\documentclass[12pt]{amsart}%
\usepackage{amsfonts}
\usepackage{mathrsfs}
\usepackage{fancyhdr}
\usepackage{comment}
\usepackage{commath}
\usepackage[a4paper, top=2.5cm, bottom=2.5cm, left=2.2cm, right=2.2cm]%
{geometry}
\usepackage{times}
\usepackage{amsmath}
\usepackage{changepage}
\usepackage{mathtools}
\usepackage{amssymb}
\usepackage{graphicx}%
\usepackage{amsthm}
\usepackage{tikz}
\usepackage{verbatim}
\usepackage{amsrefs}
\usepackage{hyperref}

\usetikzlibrary{arrows.meta}
\usetikzlibrary{positioning}
\newtheorem{theorem}{Theorem}[section]

\newtheorem{lemma}[theorem]{Lemma}

\newtheorem{definition}[theorem]{Definition}

\newtheorem{proposition}[theorem]{Proposition}
\newtheorem{remark}[theorem]{Remark}
\newtheorem{assumption}{A}

\newcommand{\N}{\mathbb{N}}

\newcommand{\R}{\mathbb{R}}

\newcommand{\B}{{\partial\Omega}}
\newcommand{\n}{\vec{n}}

\newcommand{\bb}[1]{\mathbb{#1}}

\title{Mean Field Games of Controls with Boundary Conditions \& Invariance Constraints}
\author{P.~Jameson Graber}
\thanks{The authors are grateful to be supported by National Science Foundation through NSF Grant DMS-2045027.}
	\address{J.~Graber: Baylor University, Department of Mathematics;\\
		Sid Richardson Building\\
		1410 S.~4th Street\\
		Waco, TX 76706
	}
	\email{Jameson\_Graber@baylor.edu}

\author{Kyle Rosengartner}
\address{K.~Rosengartner: Baylor University, Department of Mathematics;\\
		Sid Richardson Building\\
		1410 S.~4th Street\\
		Waco, TX 76706
	}
    \email{Kyle\_Rosengartner1@baylor.edu}
\date{December 2025}

\begin{document}

\begin{abstract}
    In a mean field game of controls, a large population of identical players seek to minimize a cost that depends on the joint distribution of the states of the players and their controls. We first consider the classes of mean field games of controls in which the value function and the distribution of player states satisfy either Dirichlet or Neumann boundary conditions. We prove that such systems are well-posed either with sufficient smallness conditions or in the case of monotone couplings. Next, we consider mean field games of controls under invariance constraints imposed on the state space. We prove the existence and uniqueness of weak solutions to our mean field game system, and then we prove higher regularity of solutions under some additional assumptions.
\end{abstract}

\maketitle

\section{Introduction}
\label{Sec: Intro}

A mean field game (MFG) is a type of differential game in which a large population of identical rational players seek to minimize a cost (or maximize a utility) that depends on the distribution of player states. The theory of mean field games was introduced independently by Lasry and Lions in \cite{lasry2007mean} and by Caines, Huang, and Malham\'e in \cite{huang2006large}. A typical mean field game can be characterized by a forward-backward system of PDE in which the optimal cost $u$ for a generic agent at time $t$ and state $x$ satisfies a backward-in-time Hamilton-Jacobi (HJ) equation and the distribution $m$ of player states satisfies a forward-in-time Fokker-Planck (FP) equation:
\begin{equation}
    \begin{cases}
        -u_t - \sum_{i,j} a_{ij}(x)\partial_{ij}^2 u + H(t,x,D_xu) = f(t,x,m), &(t,x) \in [0,T] \times \overline{\Omega} \\
        m_t - \sum_{i,j} \partial_{ij}^2(a_{ij}(x) m) - \nabla \cdot (mD_pH(t,x,D_x u)) = 0, &(t,x) \in [0,T] \times \overline{\Omega} \\
        u(T,x) = g(x,m(T)), \hspace{1cm} m(0,x) = m_0(x), &x \in \overline{\Omega}
    \end{cases}
    \label{Eq: Typical MFG}
\end{equation}
In System \eqref{Eq: Typical MFG}, the coupling is only through the distribution of players' states, i.e.~through $m$.

By contrast, in mean field games of controls (MFGCs), each player's cost depends on the joint distribution $\mu$ of states and controls (see System \eqref{Eq: MFGC}). This type of game is elsewhere referred to as an \textit{extended mean field game} (see \cites{gomes2016extended,gomes2014existence}), but the terminology ``mean field game of controls" now appears to be standard, cf.~\cite{cardaliaguet2018mean}. Compared to System \eqref{Eq: Typical MFG}, MFGCs have been much less studied in the literature. 
Perhaps the most comprehensive results on the existence of solutions to MFGCs are found in Kobeissi's 2022 papers \cites{kobeissi2022mean,kobeissi2022classical}. To derive these results, it is necessary to make a detailed study of how the Hamiltonian $H$ depends on the distribution of states and especially controls, in comparison with its dependence on the momentum variable $D_x u$, since this comparison will determine which a priori estimates are possible. A certain class of ``potential'' MFGCs has been studied in \cites{bonnans2019schauder,graber2021weak}, where the potential provides an alternative way to establish a priori estimates using the calculus of variations. The existence of solutions is now well-established on the torus $\mathbb{T}^n$ (i.e.~with periodic boundary conditions) or on $\R^n$. 

The focus of the first sections of the present article is on MFGCs with Dirichlet and Neumann boundary conditions on a smooth bounded domain.
Such problems naturally arise in applications such as economics \cites{achdou2014partial,gueant2011mean}.
Dirichlet boundary conditions model players who must leave the game on reaching a certain threshold.
This is the case, for example, in the exhaustible resource production models found in \cites{chan2015bertrand,chan2017fracking}, which have inspired a number of mathematical results for particular classes of MFGCs with Dirichlet and Neumann boundary conditions \cites{graber2018existence,graber2018variational,graber2020mean,graber2023master,camilli2025learning}.
Our purpose is to present a more general theory of MFGCs with boundary conditions.

A recent paper by Bongini and Salvarani has addressed the existence of solutions to mean field games of controls under Dirichlet boundary conditions \cite{bongini2024mean}.
The authors assumed that the set of controls is compact, and therefore the Hamiltonian is linearly bounded.
By contrast, we wish to examine the case where the Hamiltonian is coercive with respect to the momentum variable.
As for Neumann boundary conditions, some probabilistic results are given by \cite{bo2025mean} for mean field games of controls with state reflections.
Beyond this, we are unaware of any previous works on MFGC with Neumann boundary conditions. Finally, some have examined invariance conditions for MFGs (see \cite{porretta2020mean}) as well as the master equation (see \cite{zitridis2022master}). However, to our knowledge, there has been no investigation of MFGCs under invariance conditions for the state space.

The purpose of the first part of this article is to prove the well-posedness of mean field games of controls with both Dirichlet and Neumann boundary conditions.
We prove that our results for Dirichlet and Neumann boundary conditions hold under two different sets of assumptions on the Hamiltonian and/or Lagrangian.
In the first set, we carefully parametrize the growth of the Hamiltonian with respect to the distribution of controls and relate it to the growth with respect to the momentum variable, in the spirit of \cite{kobeissi2022classical}.
Uniqueness is then proved under a smallness assumption on the data.
In the second set, we impose the well-known Lasry-Lions monotonicity condition, which in the case of MFGCs is most conveniently imposed on the Lagrangian $L$ rather than the Hamiltonian $H$, under the standard assumption that $H$ is the Legendre transform of $L$ with respect to the velocity variable; cf.~\cites{kobeissi2022mean,cardaliaguet2018mean}.
This structure allows an alternative way to prove many of the a priori estimates leading to the existence of solutions, and in addition it guarantees uniqueness of solutions without any smallness assumptions.

In the last section of the article, we consider the class of MFGCs under invariance constraints, in which we impose conditions on the drift-diffusion terms such that the domain $\Omega$ is an invariant set for the controlled dynamics of the players, regardless of their controls.
In the control community, this property is sometimes referred to as the \textit{viability of the state space}. In the spirit of \cite{porretta2020mean}, most of our investigation will be done by taking solutions to an approximating MFGC on a sequence of subdomains. This will rely on the well-posedness of MFGCs under Neumann boundary conditions, and we will focus our attention in this section on the case of monotone coupling.

This type of invariance constraint is a special case of \textit{state constraints}, under which the dynamical state is forced to remain inside the domain (with probability $1$), often by putting restrictions on the class of admissible controls.
In \cite{cannarsa2018existence,cannarsa2021mean}, the authors examine the first-order mean field game system under state constraints. More recent papers (see \cite{porretta2023ergodic,sardarli2021ergodic}) have considered the second-order ergodic MFG system with state constraints (including infinite Dirichlet boundary conditions), which relied on the aysmptotic behavior of $D_xu$ near the boundary, established in \cite{alessio2006asymptotic}. In \cite{lasry1989nonlinear}, Lasry and Lions establish properties of solutions to elliptic equations with state constraints in both the ergotic and non-ergotic cases.

Our main contribution is as follows.
We provide a systematic, comprehensive set of results on existence and uniqueness to MFGCs with both Dirichlet and Neumann boundary conditions.
We synthesize many of the ideas found in Kobeissi's papers \cites{kobeissi2022classical,kobeissi2022mean} so as to provide one self-contained treatment of several classes of coupling, both monotone and non-monotone.
More importantly, we provide a priori estimates on classical solutions that apply even when the boundary conditions create additional difficulties.
For example, to prove an a priori bound on the gradient of the value function is considerably more technical, since to apply a Bernstein type argument one has to analyze its behavior near the bounday (see Section \ref{sec: Gradient Bounds}).
As for the distribution of states and controls, Dirichlet boundary conditions create the added difficulty of mass absorption, so that the standard Wasserstein metric is no longer the appropriate tool to analyze the behavior of the distribution of states.
We handle this differently from \cite{bongini2024mean}, preferring to use a metric very much akin to the Wasserstein metric (cf.~\cite{graber2023master}). Additionally, we combine these results with the methods used in \cite{porretta2020mean} to give a number of results on the existence, uniqueness, and regularity solutions to MFGCs under invariance constraints on the state space.

The remainder of this manuscript is organized as follows.
In the rest of this introduction, we introduce some notation and useful preliminaries, followed by a presentation of the PDE systems that we study.
In Section \ref{sec: Problem 1}, we introduce the first set of assumptions and prove existence of solutions to \eqref{Eq: MFGC} under Dirichlet and Neumann boundary conditions with non-monotone coupling, under carefully parametrized smallness assumptions.
Then in Section \ref{sec: Problem 2}, we introduce the second set of assumptions, in particular the Lasry-Lions monotonicity condition, and prove the existence of solutions to this system under such couplings.
In Section \ref{Sec: Uniqueness}, we prove the uniqueness of solutions under both monotone and non-monotone coupling.
Finally, in Section \ref{sec: Problem 3}, we introduce the assumptions that we will use in our analysis of MFGCs under invariance constraints; we prove the existence and uniqueness of solutions, as well as increased regularity of solutions under additional assumptions.

\subsection{Notation and Preliminaries}
\label{Sec: Notation}

We begin by recalling a few definitions from \cite{porretta2020mean}.
\begin{definition} \label{Def: Compact Domain of Class C2}
    We say $K \subseteq \R^n$ is a compact domain of class $C^2$ if $K$ is a compact connected set and there exists $M \in \N$ such that for $1 \leq i \leq M$, there exist $x_i \in \partial K$, $r_i > 0$, $\phi_i: B_{r_i}(x_i) \rightarrow \R$ so that
    \begin{enumerate}
        \item $\partial K \subseteq \bigcup_{i=1}^M B_{r_i}(x_i)$
        \item $\partial K \cap B_{r_i}(x_i) = \{x \in B_{r_i}(x_i) : \phi_i(x) = 0\}$
        \item $\phi_i \in C^2(B_{r_i}(x_i))$ with $D_{xx}^2\phi_i$ bounded.
    \end{enumerate}
\end{definition}

\begin{definition} \label{Def: Subdomains}
    For every $\varepsilon > 0$, we will denote by $\Gamma_\varepsilon$ and $\Omega_\varepsilon$ the sets
    $$\Gamma_\varepsilon \coloneqq \{x \in \R^n : d(x,\Omega) \leq \varepsilon\} \hspace{1cm} \Omega_\varepsilon \coloneq \Omega \setminus \Gamma_\varepsilon.$$
    Furthermore, we will denote by $d_\Omega$ a $C^2(\overline{\Omega})$ function such that there exists $\varepsilon_0 > 0$ with $d_\Omega(x) = d(x,\Omega)$ in $\Gamma_{\varepsilon_0}$ (see \cite{cannarsa2010invariant,delfour1994shape}). When there is no ambiguity, we will merely write $d$ for $d_\Omega$.
\end{definition}

Next, we define the following spaces of measures and the metrics we impose on them, which will vary depending on the boundary conditions.
Given $R > 0$ and sets $A,B \subseteq \R^n$, define $\mathfrak{M}(A)$ to be the set of Borel measures $\rho$ on $A$ with $\rho(A) \leq 1$ and define $\mathfrak{M}_{\infty,R}(A \times B)$ to be the set of Borel measures $\rho \in \mathfrak{M}(A \times B)$ with $\text{supp}\rho \subseteq \{(a,b) \in A \times B : |b| \leq R\}$.
We endow these spaces with a metric, depending on the boundary conditions.
In the Neumann case, as the ``total mass" $\int_\Omega m(t) dx$ remains constant, we can restrict to the space $\mathcal{P} \subset \mathfrak{M}$ consisting of probability measures, which is endowed with the Wasserstein metric
\begin{equation}
	d^*(\mu,\nu) = W_1(\mu,\nu) = \inf\left\{\int \abs{x-y} d\pi(x,y) : \pi \in \Pi(\mu,\nu)\right\}
\end{equation}
where $\Pi(\mu,\nu)$ denotes the set of all couplings between $\mu$ and $\nu$.
By Kantorovitch duality, we also have the characterization
\begin{equation}
	W_1(\mu,\nu) = \sup\left\{\int \varphi d(\mu-\nu) : \|\varphi\|_{Lip}  \leq 1\right\}.
\end{equation}
See \cite{villani2021topics} for more details.
We will also denote by $\mathcal{P}_{\infty,R}(A \times B)$ the set of all probability measures in $\mathfrak{M}_{\infty,R}(A \times B)$.

To deal with the ``mass escape" that occurs in the Dirichlet case and the variation of mass in the sequence of solutions to approximating problems found in Section \ref{sec: Problem 3}, we will endow $\mathfrak{M}$ with the metric
$$d^*(\mu,\nu) = \sup\left\{\int \varphi d(\mu-\nu) : \|\varphi\|_{Lip}  \leq 1, \|\varphi\|_\infty \leq 1\right\}.$$

We will also define the quantities
\begin{equation}
    \Lambda_q(\mu) \coloneqq \left(\int_{\Omega \times \R^n} |\alpha|^q d\mu(x,\alpha)\right)^\frac{1}{q}
    \label{Eq: Lambda-q}
\end{equation}
for $1 \leq q < \infty$ and
\begin{equation}
    \Lambda_\infty(\mu) \coloneqq \sup\{|\alpha| : (x,\alpha) \in \text{supp} \mu\}
    \label{Eq: Lambda-infinity}
\end{equation}
as in \cite{kobeissi2022classical,kobeissi2022mean}. These will allow us to quantify the dependence of $H$ on $\mu$.

As for function spaces, for non-negative integers $j,k$ we denote by $C^{j,k}$ the space of all functions $u(t,x)$ on $Q$ that are $j$ times continuously differentiable with respect to $t$ and $k$ times continuously differentiable with respect to $x$.
For a fraction $\alpha \in (0,1)$ we denote by $C^{1+\alpha/2,2+\alpha}$ the usual parabolic H\"older space as introduced for instance in \cite{ladyzhenskaia1968linear}.
As for Sobolev spaces, we use similar notation as in \cite{ladyzhenskaia1968linear}, in particular $W^{1,2}_p$ denotes the space of all functions with weak derivatives up to order 1 in time and 2 in space, whose weak derivatives are all $L^p$ summable.

Finally, we recall the well-known Leray-Schauder fixed point theorem by which we will prove existence of solutions.

\begin{theorem}[Leray-Schauder]
    Let $X$ be a Banach space and let $T: X \times [0,1] \rightarrow X$ be a continuous and compact mapping. Assume there exist $x_0 \in X$ and $C > 0$ so that $T(x,0) = x_0$ for all $x \in X$ and $\|x\|_X < C$ for all $(x,\tau) \in X \times [0,1]$ such that $T(x,\tau) = x$. Then there exists $x \in X$ such that $T(x,1) = x$.
\end{theorem}

Aside from these preliminaries, we specify that the constant $C$ appearing in many results denotes a generic constant that depends only on the data, specifically on the constants found in the Assumptions (that is, in Section \ref{sec:assms1} or \ref{sec:assms2}).

\subsection{The Systems of PDE}
\label{Sec: Our MFGC}

Let $\Omega \subseteq \R^n$ be a bounded, convex (we need not assume convexity in the Dirichlet case), $C^{2+\beta}$ domain for some $\beta \in (0,1)$. In Sections \ref{sec: Problem 1} and \ref{sec: Problem 2}, we will consider the system
\begin{equation}
    \begin{cases}
        -u_t - \nu\Delta u + H(t,x,D_xu,\mu) = f(t,x,m), &(t,x) \in Q \\
        m_t - \nu\Delta m - \nabla \cdot (mD_pH(t,x,D_x u,\mu)) = 0, &(t,x) \in Q \\
        \mu = (I,-D_pH(t,\cdot,D_xu,\mu)) \# m, &t \in [0,T] \\
        u(T,x) = g(x,m(T)), \hspace{1cm} m(0,x) = m_0(x), &x \in \overline{\Omega}
    \end{cases}
    \label{Eq: MFGC}
\end{equation}
paired with either Dirichlet
\begin{equation}
    u = m = 0, \hspace{1cm} (t,x) \in \Sigma \tag{\ref{Eq: MFGC}d}
    \label{Eq: Dirchlet}
\end{equation}
or Neumann
\begin{equation}
    \frac{\partial u}{\partial\n} = \nu\frac{\partial m}{\partial\n} + mD_pH(t,x,D_xu,\mu) \cdot \n = 0, \hspace{1cm} (t,x) \in \Sigma \tag{\ref{Eq: MFGC}n}
    \label{Eq: Neumann}
\end{equation}
boundary conditions, where $Q \coloneqq [0,T] \times \overline{\Omega}$ and $\Sigma \coloneqq [0,T] \times \B$ for some $T > 0$.

For the purpose of variation, we will investigate classical solutions in Section \ref{sec: Problem 1} and strong solutions in Section \ref{sec: Problem 2}. We define these respectively as follows.

\begin{definition} \label{Def: Classical Solution}
    We will say $(u,m,\mu)$ is a classical solution to \eqref{Eq: MFGC}-\eqref{Eq: Dirchlet} (resp., \eqref{Eq: MFGC}-\eqref{Eq: Neumann}) if
    \begin{enumerate}
        \item $u \in C^{1,2}(Q)$ is a classical solution to the Hamilton-Jacobi equation;
        \item $m \in \{\rho \in L^2(0,T;H_0^1(\Omega)) : \rho_t \in L^2(0,T;H^{-1}(\Omega))\}$ (resp.~$m \in \{\rho \in L^2(0,T;H^1(\Omega)) : \rho_t \in L^2(0,T;H^1(\Omega)^*)\}$) satisfies
        \begin{equation} \label{Eq: Weak FP}
            \int_\Omega (m_t\varphi + (\nu D_xm + mD_pH) \cdot D_x\varphi)dx = 0
        \end{equation}
        for all $\varphi \in H^1_0(\Omega)$ (resp., $\varphi \in H^1(\Omega)$) and a.e. $0 \leq t \leq T$;
        \item $\mu \in C^0(0,T;\mathfrak{M}(\overline{\Omega} \times \R^n))$ (resp.~$\mu \in C^0(0,T;\mathcal{P}(\overline{\Omega} \times \R^n))$) satisfies the fixed-point relation $\mu = (I,-D_pH(t,\cdot,D_xu,\mu)) \# m$ at every $t \in [0,T]$;
        \item $(u,m)$ satisfies $u(T,\cdot) = g(\cdot,m(T))$, $m(0,\cdot) = m_0$, and $u|_\Sigma = 0$ (resp., $\frac{\partial u}{\partial \n}|_\Sigma = 0$) pointwise.
    \end{enumerate}
\end{definition}

\begin{definition} \label{Def: Strong Solution}
    We will say $(u,m,\mu)$ is a strong solution to \eqref{Eq: MFGC}-\eqref{Eq: Dirchlet} (resp., \eqref{Eq: MFGC}-\eqref{Eq: Neumann}) if
    \begin{enumerate}
        \item $u \in C^{0,1}(Q) \cap W^{1,2}_2(Q)$ is a strong solution to the HJ equation, i.e.,
        $$-u_t - \nu\Delta u + H(t,x,D_xu,\mu) = f(t,x,m)$$
        for a.e. $(t,x) \in Q$;
        \item $m \in \{\rho \in L^2(0,T;H_0^1(\Omega)) : \rho_t \in L^2(0,T;H^{-1}(\Omega))\}$ (resp.~$m \in \{\rho \in L^2(0,T;H^1(\Omega)) : \rho_t \in L^2(0,T;H^1(\Omega)^*)\}$) satisfies \eqref{Eq: Weak FP}
        for all $\varphi \in H^1_0(\Omega)$ (resp., $\varphi \in H^1(\Omega)$) and a.e. $0 \leq t \leq T$;
        \item $\mu \in C^0(0,T;\mathfrak{M}(\overline{\Omega} \times \R^n))$ satisfies $\mu = (I,-D_pH(t,\cdot,D_xu,\mu)) \# m$ at every $t \in [0,T]$;
        \item $(u,m)$ satisfies $u(T,\cdot) = g(\cdot,m(T))$, $m(0,\cdot) = m_0$, and $u|_\Sigma = 0$ (resp., $\frac{\partial u}{\partial \n}|_\Sigma = 0$) pointwise.
    \end{enumerate}
\end{definition}

It will be useful in later sections to notice the following stochastic interpretations of our boundary value problems.

\begin{remark} \label{Rmk: Stochastic Interpretation}
    We note that $X(t) \sim m(t)$ and $(X(t),-D_pH(t,X(t),D_xu(t,X(t)),\mu(t))) \sim \mu(t)$, where
    \begin{enumerate}
        \item In the Dirichlet case,
        $$X(t) = x_0 - \int_0^{t \land \tau} D_pH(s,X(s),D_xu(s,X(s)),\mu(s)) ds + \sqrt{2\nu}\int_0^{t \land \tau} dB(s),$$
        where $B(t)$ denotes $n$-dimensional Brownian motion and $\tau$ is the stopping time
        $$\tau \coloneqq \inf\{t \in [0,T] : X(t) \in \B \text{ or } t = T\}$$
        (See \cite{bongini2024mean}).
        \item In the Neumann case,
        $$
        \begin{cases}
            dX(t) = -D_pH(t,X(t),D_xu(t,X(t)),\mu(t))dt + \sqrt{2\nu}dB(t) - \int_{\partial \Omega} \n(X(t)) dL(t) \\
            X(0) = x_0
        \end{cases}
        $$
        where $L(t)$ is the local time of $X(t)$ and $\n(x)$ is the unit outward normal vector.
        For convex $\Omega$, \cite{tanaka1979stochastic} shows that this equation can be solved. See \cite{bo2025mean} for an analysis of reflection boundary conditions in the case of MFGCs specifically.
    \end{enumerate}
\end{remark}

In Section \ref{sec: Problem 3}, we will show that the system
\begin{equation}
    \begin{cases}
        -u_t - \sum_{i,j} a_{ij}(x)\partial_{ij}^2 u + H(t,x,D_xu,\mu) = f(t,x,m), &(t,x) \in Q \\
        m_t - \sum_{i,j} \partial_{ij}^2(a_{ij}(x) m) - \nabla \cdot (mD_pH(t,x,D_x u,\mu)) = 0, &(t,x) \in Q \\
        \mu = (I,-D_pH(t,\cdot,D_xu,\mu)) \# m, &t \in [0,T] \\
        u(T,x) = g(x,m(T)), \hspace{1cm} m(0,x) = m_0(x), &x \in \overline{\Omega}
    \end{cases}
    \label{Eq: MFGC 2}
\end{equation}
is well-posed in a weak sense without need for boundary conditions, given the assumption that the diffusion coefficient $a$ and the Hamiltonian $H$ satisfy the following invariance constraint:
\begin{equation} \tag{\ref{Eq: MFGC 2}*}
    tr(a(x)D^2d(x)) - D_pH(t,x,p,\mu) \cdot Dd(x) \geq \frac{a(x)Dd(x) \cdot Dd(x)}{d(x)} - Cd(x)
    \label{Eq: Invariance Condition}
\end{equation}
for all $t \in [0,T]$ and $p \in \R^n$, for some $C > 0$, for $x$ in some neighborhood of $\B$. To this end, we define weak solutions as in \cite{porretta2020mean}.

\begin{definition} \label{Def: Solutions to HJ}
    We will say $u$ is a weak solution of
    \begin{equation} \label{Eq: HJ}
        \begin{cases}
            -u_t - \sum_{i,j} a_{ij}(x)\partial_{ij}^2 u + H(t,x,D_xu) = 0, &(t,x) \in Q \\
                u(T,x) = g(x), &x \in \Omega
        \end{cases}
    \end{equation}
    provided
    \begin{enumerate}
        \item $u \in L^\infty(Q)$;
        \item $u \in L^2(0,T;W^{1,2}(K))$ for each $K \subset\subset \Omega$
        \item For every $\varphi \in C^\infty((0,T] \times \Omega)$,
        $u$ satisfies
        $$\int_0^T\int_\Omega \left[u\varphi_t + aD_xu \cdot D_x\varphi + (\Tilde{b}D_xu + H)\varphi\right] dxdt = \int_\Omega g\varphi(T) dx$$
        where $\Tilde{b}_j(x) = \sum_{i=1}^n \partial_i a_{ij}(x)$ for $j = 1,\dots,n$.
    \end{enumerate}
\end{definition}

\begin{definition} \label{Def: Solutions to FP}
    Given a locally bounded vector field $b: [0,T] \times \Omega \rightarrow \R^n$, we will say $m \in L^1([0,T] \times \Omega)$ is a weak solution of \begin{equation}
        \begin{cases}
            m_t - \sum_{i,j} \partial_{ij}^2 (a_{ij}(x)m) - \nabla \cdot (mb) = 0, &(t,x) \in (0,T) \times \Omega \\
            m(0,x) = m_0(x), &x \in \Omega
        \end{cases}
        \label{Eq: FP}
    \end{equation}
    provided
    \begin{enumerate}
        \item $m \in C^0(0,T; L^1(\Omega))$ with $m \geq 0$;
        \item For every $\varphi \in C(0,T;L^1(\Omega)) \cap L^\infty(Q)$ satisfying
        $$\begin{cases}
            -\varphi_t - \sum_{i,j} a_{ij}\partial_{ij}^2\varphi + b \cdot D_x\varphi \in L^\infty(Q) \\
            \varphi(T) = 0
        \end{cases}$$
        in the sense of distributions, we have
        $$\int_0^T\int_\Omega m\left(-\varphi_t - \sum_{i,j} a_{ij}\partial_{ij}^2\varphi + b \cdot D_x\varphi\right)dxdt = \int_\Omega m_0\varphi(0)dx.$$
    \end{enumerate}
\end{definition}

\begin{definition} \label{Def: Weak Solution}
    Given $R > 0$, we will consider $(u,m,\mu) \in C^{0,1}([0,T] \times \Omega) \times C^0(0,T; L^1(\Omega))_+ \times C^0(0,T;\mathfrak{M}_{\infty,R}(\Omega \times \R^n))$ a weak solution to \eqref{Eq: MFGC 2} if
    \begin{enumerate}
        \item $u$ is a weak solution to the first line of \eqref{Eq: MFGC 2} in sense of Definition \ref{Def: Solutions to HJ};
        \item $m$ is a weak solution to the second line of \eqref{Eq: MFGC 2} in sense of Definition \ref{Def: Solutions to FP};
        \item $\mu$ satisfies the third equation of \eqref{Eq: MFGC 2} at every $t \in [0,T]$;
        \item $(u,m)$ satisfy the last line of \eqref{Eq: MFGC 2} pointwise.
    \end{enumerate}
\end{definition}

\section{Dirichlet \& Neumann Boundary Conditions with Non-Monotone Coupling} \label{sec: Problem 1}

For the first problem, we will make the majority of our assumptions on the Hamiltonian, and we will prove existence of classical solutions using an approach similar to that of \cite{kobeissi2022classical}.

With the aim of applying the Leray-Schauder fixed point theorem, we will consider the following parametrized system:
\begin{equation}
    \begin{cases}
        -u_t - \nu\Delta u + \lambda H(t,x,D_xu,\mu) = 0, &(t,x) \in Q \\
        m_t - \nu\Delta m - \lambda\nabla \cdot (mD_pH(t,x,D_x u,\mu)) = 0, &(t,x) \in Q \\
        \mu = (I,-\lambda D_pH(t,\cdot,D_xu,\mu)) \# m, &t \in [0,T] \\
        u(T,x) = \lambda g(x,m(T)), \hspace{1cm} m(0,x) = \lambda m_0(x), &x \in \overline{\Omega}
    \end{cases}
    \label{Eq: Parametrized MFGC 1}
\end{equation}
paired with either
\begin{equation}
    u = m = 0, \hspace{1cm} (t,x) \in \Sigma \tag{\ref{Eq: Parametrized MFGC 1}d}
    \label{Eq: Parametrized Dirchlet 1}
\end{equation}
or
\begin{equation}
    \frac{\partial u}{\partial\n} = \nu\frac{\partial m}{\partial\n} + \lambda mD_pH(t,x,D_xu,\mu) \cdot \n = 0, \hspace{1cm} (t,x) \in \Sigma \tag{\ref{Eq: Parametrized MFGC 1}n}
    \label{Eq: Parametrized Neumann 1}
\end{equation}
where $\lambda \in [0,1]$ is given.

\subsection{Assumptions} \label{sec:assms1}

In the case of non-monotone couplings, we will make the following assumptions.
The constants $C_0,\lambda_0,\lambda_1,\lambda_2,q,q_0$ and the functions $\xi_1,\xi_2,\xi_3$ listed below are fixed independent of the data.
\begin{assumption} \label{A: Hamiltonian}
    For every $R > 0$, the Hamiltonian $H: [0,T] \times \overline{\Omega} \times \R^n \times \mathfrak{M}_{\infty,R}(\Omega \times \R^n) \rightarrow \R$ is differentiable with respect to $(x,p)$, strictly convex with respect to $p$, satisfies the coercivity condition
    $$\lim_{|p| \rightarrow \infty} \frac{H(t,x,p,\mu)}{|p|} = +\infty,$$
    and $H$ and its derivatives are continuous with respect to $(t,x,p,\mu)$.
    Moreover, the Lagrangian $L: [0,T] \times \overline{\Omega} \times \R^n \times \mathfrak{M}_{\infty,R}(\Omega \times \R^n) \rightarrow \R$ is $C^1$ and strictly convex with respect to $\alpha$.
\end{assumption}
\begin{assumption} \label{A: Initial & Terminal Conditions}
    We take $m_0 \in C^\beta(\overline{\Omega})_+$
    with $\int_{\Omega} m_0 dx = 1$. For all $m \in \mathfrak{M}(\overline{\Omega})$, we have $g(\cdot,m) \in C^{3 + \beta}(\overline{\Omega})$ with $g|_\B = 0$ (resp., $\frac{\partial g}{\partial\n}|_\B = 0$) and $\|g(\cdot,m)\|_{C^{3+\beta}(\overline{\Omega})} \leq C_0$ for some constant $C_0$.
    Furthermore, $m \mapsto g(\cdot,m)$ is continuous from $L^2(\Omega)$ into $L^2(\Omega)$.
\end{assumption}
\begin{assumption} \label{A: f is zero}
    In this problem, we will assume $f \equiv 0$. (If $f$ is sufficiently smooth such that $f$ and its derivatives are bounded, we need only replace $H$ with $H - f$.
    Then there is no loss of generality in assuming $f \equiv 0$.)
\end{assumption} 
\begin{assumption} \label{A: Bound for DpH}
    $|D_pH(t,x,p,\mu)| \leq C_0(1 + |p|^{q-1}) + \lambda_0\Lambda_{q_0}(\mu)$ for some constants $q \in (1,\infty)$, $q_0 \leq q' \coloneqq q/(q-1)$, and $\lambda_0 \in (0,1)$; in the Dirichlet case we require $q \leq 2$.
\end{assumption}
\begin{assumption} \label{A: Convexity assumption for H}
    $D_pH(t,x,p,\mu) \cdot p - H(t,x,p,\mu) \geq C_0^{-1}(|p|^q - \lambda_1\Lambda_{q_0}(\mu)^{q'}) - C_0$.
\end{assumption}
\begin{assumption} \label{A: Bounding H}
    $|H(t,x,0,\mu)| \leq \lambda_2\Lambda_{q_0}(\mu)^{q'} + C_0$, where $\lambda_1$ and $\lambda_2$ satisfy $0 \leq \lambda_1 < \frac{(1-\lambda_0)^{q'}}{C_0^{q'}} - C_0\lambda_2$.
\end{assumption}
\begin{assumption} \label{A: Bound for DxH}
    $|D_xH(t,x,p,\mu)| \leq C_0(1 + |p|^q + \Lambda_{q_0}(\mu)^{q'})$.
\end{assumption}
\begin{assumption} \label{A: DpH is Lipschitz in alpha}
    For every $R > 0$, $p \in \R^n$, and $\mu_1,\mu_2 \in \mathfrak{M}_{\infty,R}(\Omega \times \R^n)$ with $\mu_i = (I,\alpha^{\mu_i}) \# \rho$, we have
    $$|D_pH(t,x,p,\mu_1) - D_pH(t,x,p,\mu_2)| \leq L_1\|\alpha^{\mu_1} - \alpha^{\mu_2}\|_{L^{q_0}(\rho)}$$
    for some $L_1 \in (0,1)$.
\end{assumption}
\begin{assumption} \label{A: DpH is Lipschitz in m}
    For every $R > 0$, $p \in \R^n$, and $\mu_1,\mu_2 \in \mathfrak{M}_{\infty,R}(\Omega \times \R^n)$ with $\mu_i = (I,\alpha) \# \rho_i$, we have
    $$|D_pH(t,x,p,\mu_1) - D_pH(t,x,p,\mu_2)| \leq d^*(\rho_1,\rho_2)^\beta\xi_1(\abs{p},\|\alpha\|_\infty)$$
    for some continuous function $\xi_1: \R^2 \rightarrow [0,\infty)$.
\end{assumption}
\begin{assumption} \label{A: Regularity of DpH}
    For every $R > 0$, $p \in \R^n$, and $\mu \in \mathfrak{M}_{\infty,R}(\Omega \times \R^n)$, we have
    $$
    \begin{aligned}
        &|D_pH(t,x,p_1,\mu) - D_pH(s,y,p_2,\mu)| \\
        &\leq (|p_1 - p_2|^\beta + |x-y|^\beta + |t-s|^{\beta/2})\xi_2(\abs{p_1},\abs{p_2},\Lambda_\infty(\mu))
    \end{aligned}
    $$
    for some continuous function $\xi_2: \R^3 \rightarrow [0,\infty)$.
\end{assumption}
\begin{assumption} \label{A: Regularity of H}
    For every $R > 0$, $p \in \R^n$, and $\mu_1,\mu_2 \in \mathfrak{M}_{\infty,R}(\Omega \times \R^n)$, we have
    $$|H(t,x,p,\mu_1) - H(s,x,p,\mu_2)| \leq (|t-s|^{\beta/2} + d^*(\mu_1,\mu_2)^\beta)\xi_2(\abs{p},\Lambda_\infty(\mu_1),\Lambda_\infty(\mu_2)).$$
\end{assumption}
\begin{assumption} \label{A: Regularity of DxH}
    For every $R > 0$, $p_1,p_2 \in \R^n$, and $\mu_1,\mu_2 \in \mathfrak{M}_{\infty,R}(\Omega \times \R^n)$, we have
    $$
    \begin{aligned}
        &|D_xH(t,x,p_1,\mu) - D_xH(s,y,p_2,\mu)| \\
        &\leq (|x-y|^\beta + |t-s|^{\beta/2} + |p_1 - p_2|^\beta + d^*(\mu_1,\mu_2)^\beta)\xi_3(\abs{p_1},\abs{p_2},\Lambda_\infty(\mu_1),\Lambda_\infty(\mu_2))
    \end{aligned}
    $$
    for some continuous function $\xi_3: \R^4 \rightarrow [0,\infty)$.
\end{assumption}

For motivating examples, one could pair the examples in \cite[Section 6]{kobeissi2022classical} with relevant boundary conditions. For instance, one could consider the linear-demand exhaustible resource model with non-positively correlated resources paired with Dirichlet boundary conditions, which would model a situation in which players must leave the game when their production capacities reach the boundary of a specified region. Alternatively, one consider models of crowd dynamics or flocking birds paired with Neumann boundary conditions, which would model a situation in which members are reflected off of the boundary and must therefore remain inside the enclosed region. It is straightforward to check that our assumptions are satisfied by these examples.

\subsection{Fixed-Point Relation in $\mu$}
\label{Sec: Fixed-Point Relation}
In this section, we address the well-posedness of the third equation in System \eqref{Eq: MFGC}, which we regard as a fixed-point relation for the measure $\mu$.

\begin{lemma} \label{Lem: Existence & Bound for mu 1}
    Assume A\ref{A: Hamiltonian}, A\ref{A: Bound for DpH}, and A\ref{A: DpH is Lipschitz in alpha} hold. Given $t \in [0,T]$, $\lambda \in [0,1]$, $p \in C^0(\Omega;\R^n)$, and $m \in \mathfrak{M}(\overline{\Omega})$, we have the following:
    \begin{enumerate}
        \item[1)] There exists a unique $\mu \in \mathfrak{M}(\overline{\Omega} \times \R^n)$ such that
        \begin{equation}
            \mu = (I,-\lambda D_pH(t,\cdot,p(\cdot),\mu)) \# m.
            \label{Eq: mu 1}
        \end{equation}
        \item[2)] For any $\Tilde{q} \in [1,\infty]$, we have
        \begin{equation}
            \Lambda_{\Tilde{q}}(\mu) \leq \frac{\lambda C_0}{1-\lambda\lambda_0}\left(1 + \||p|^{q-1}\|_{L^{\max\{q_0,\Tilde{q}\}}(m)}\right).
            \label{Eq: Bound for Lambda-q}
        \end{equation}
    \end{enumerate}
    Combining \eqref{Eq: Bound for Lambda-q} and A\ref{A: Bound for DpH} gives $R > 0$ (depending on $p,m$) such that $\mu \in \mathfrak{M}_{\infty,R}(\Omega \times \R^n)$.
    Moreover, if $m \in \mathcal{P}(\overline{\Omega})$ then $\mu \in \mathcal{P}(\overline{\Omega} \times \bb{R}^n)$.
\end{lemma}

\begin{remark}
    Note that since we have $\int_\Omega m dx \leq 1$, Jensen's inequality gives that
    $$\Lambda_{q_1}(\mu_\alpha) \leq \Lambda_{q_2}(\mu_\alpha)$$
    for $1 \leq q_1 \leq q_2 \leq \infty$, where $\mu_\alpha = (I,\alpha) \# m$ for some $\alpha$.
\end{remark}

\begin{proof}[Proof of Lemma \ref{Lem: Existence & Bound for mu 1}]
    1) Define $\Phi_{(p,m)}^t: C^0(\Omega;\R^n) \rightarrow C^0(\Omega;\R^n)$ given by
    $$(\Phi_{(p,m)}^t(\alpha))(x) = -\lambda D_pH(t,x,p(x),(I,\alpha) \# m).$$
    Then
    $$
    \begin{aligned}
        \|\Phi_{(p,m)}^t(\alpha_1) - \Phi_{(p,m)}^t(\alpha_2)\|_\infty &= \lambda\sup_{x \in \Omega} |D_pH(t,x,p(x),(I,\alpha_1) \# m) - D_pH(t,x,p(x),(I,\alpha_2) \# m)| \\
        &\leq L_1\|\alpha_1 - \alpha_2\|_{L^{q_0}(m)} \\
        &\leq L_1\|\alpha_1 - \alpha_2\|_\infty.
    \end{aligned}
    $$
    By the Banach fixed point theorem, since $L_1 < 1$, $\Phi_{(p,m)}^t$ has a unique fixed point $\alpha$. In particular, $\mu \coloneqq (I,\alpha) \# m$ uniquely solves \eqref{Eq: mu 1}.
    
    2) If $\Tilde{q} \geq q_0$,
    $$
    \begin{aligned}
        \Lambda_{\Tilde{q}}(\mu) &= \lambda\left\|D_pH(t,x,p,\mu)\right\|_{L^{\Tilde{q}}(m)} \\
        &\leq \lambda\left\|C_0(1 + |p|^{q-1}) + \lambda_0\Lambda_{q_0}(\mu)\right\|_{L^{\Tilde{q}}(m)} \\
        &\leq \lambda C_0(1 + \left\||p|^{q-1}\right\|_{L^{\Tilde{q}}(m)}) + \lambda\lambda_0\Lambda_{\Tilde{q}}(\mu)
    \end{aligned}
    $$
    and so
    $$\Lambda_{\Tilde{q}}(\mu) \leq \frac{\lambda C_0}{1 - \lambda\lambda_0}\left(1 + \left\||p|^{q-1}\right\|_{L^{\Tilde{q}}(m)}\right)$$
    If $\Tilde{q} \leq q_0$, we use a similar argument to get
    $$\Lambda_{\Tilde{q}}(\mu) \leq \Lambda_{q_0}(\mu) \leq \frac{\lambda C_0}{1 - \lambda\lambda_0}\left(1 + \left\||p|^{q-1}\right\|_{L^{q_0}(m)}\right).$$
\end{proof}

\begin{lemma} \label{Lem: Continuity of mu 1}
    Assume A\ref{A: Hamiltonian}-\ref{A: DpH is Lipschitz in alpha} hold. Let $(t_k,\lambda_k,p_k,m_k)_{n \in \N}$ be a sequence in $[0,T] \times [0,1] \times C^0(\overline{\Omega},\R^n) \times L^1(\Omega)$ such that
    \begin{itemize}
        \item $t_k \rightarrow t$ in $[0,T]$ and $\lambda_k \rightarrow \lambda$ in $[0,1]$;
        \item $p_k \rightarrow p$ in $C^0(\overline{\Omega};\R^n)$;
        \item $\|m_k\|_{L^1(\Omega)} \leq 1$ and $m_k \rightarrow m$ in $L^1(\Omega)$.
    \end{itemize}
    Then $\mu^k \rightarrow \mu$ in $\mathfrak{M}_{\infty,R}(\Omega \times \R^n)$, where $\mu^k$ and $\mu$ are the fixed points associated to $(t_k,\lambda_k,p_k,m_k)$ and $(t,\lambda,p,m)$, respectively.
\end{lemma}

\begin{proof}
    Since $(p_k)_{k \in \N}$ is bounded in $C^0(\overline{\Omega};\R^n)$, by (\ref{Eq: Bound for Lambda-q}),
    $$\Lambda_\infty(\mu^{k}) \leq \frac{C_0}{1 - \lambda_0}(1+\|p_k\|_{L^\infty(m_k)}^{q-1}) \leq C$$
    and so $\text{supp}\mu^{k} \subseteq \{(x,\alpha) : |\alpha| \leq C\}$. Thus, the $\mu^{k}$ are uniformly compactly supported. Hence, $(\mu^{k})_{k \in \N}$ is compact in $\mathfrak{M}_{\infty,R}(\Omega \times \R^n)$. Take a subsequence $(\mu^{k_j})_{j \in \N}$ converging to some $\Tilde{\mu}$. In the spirit of \cite[Special case 6.16]{villani2008optimal}, we have
    $$
    \begin{aligned}
        &d^*(\mu^{k_j},(I,-\lambda D_pH(t,x,p,\Tilde{\mu})) \# m) \\
        &\leq \int_{\Omega} |\lambda_{k_j}D_pH(t_{k_j},x,p_{k_j},\mu^{k_j}) - \lambda D_pH(t,x,p,\Tilde{\mu})| dm \\
        &\qquad + \max\{1,\lambda_{k_j}osc(D_pH(t_{k_j},\cdot,p_{k_j}(\cdot),\mu^{k_j}))\}\int_{\Omega} |m - m_{k_j}| dx \\
        &\rightarrow 0.
    \end{aligned}
    $$
    Therefore, $\Tilde{\mu} = (I,-\lambda D_pH(t,\cdot,p(\cdot),\Tilde{\mu})) \# m$.
    By uniqueness of the fixed point, this shows that $\mu = \Tilde{\mu}$. Thus, the map $(t,\lambda,p,m) \mapsto \mu$ is continuous.
\end{proof}

\begin{remark} \label{rmk: uniform convergence of mu 1}
    This shows that for any $(\lambda,u,m) \in [0,1] \times C^{0,1}(Q) \times C^0(0,T;L^1(\Omega))$, the map $t \mapsto \mu(t) = (I,-\lambda D_pH(t,\cdot,D_xu(t,\cdot),\mu(t)) \# m(t)$ is continuous. Furthermore, if $(\lambda_k,p_k,m_k)_{k \in \N}$ converges to $(\lambda,D_xu,m)$ in $[0,1] \times C^0(Q;\R^n) \times C^0(0,T;L^1(\Omega))$, then $\mu_k(t)$ defined by $$\mu_k(t) \coloneqq (I,-\lambda_k D_pH(t,\cdot,p_k(t,\cdot),\mu_k(t))) \# m_k(t)$$ converges to $\mu$ uniformly on $[0,T]$.
\end{remark}

\subsection{Estimates for $u$}
\label{Sec: Estimates for u}

The main purpose of this section is to establish a priori bounds on $u$ for solutions of System \eqref{Eq: MFGC}. This will require some a priori estimates on $\mu$ as well, which follow from so-called ``energy estimates" that are now standard in mean field games.

\begin{proposition}
    Assume A\ref{A: Hamiltonian}-\ref{A: f is zero} hold. Fix $\lambda \in [0,1]$ and let $(u,m,\mu)$ be a classical solution of \eqref{Eq: Parametrized MFGC 1}. Then
    $$-\max_x u(T,x)_- - \lambda\int_t^T \|H(t,\cdot,0,\mu(s))\|_\infty ds \leq u(t,x) \leq \max_x \, u(T,x)_+ + \lambda\int_t^T \|H(t,\cdot,0,\mu(s))\|_\infty ds$$
    \label{Prop: Maximum Principle}
\end{proposition}

\begin{proof}
    Note that for $w(t,x) = u(t,x) - \lambda\int_t^T \|H(t,\cdot,0,\mu(s))\|_\infty ds$, we have
    $$
    \begin{aligned}
        \frac{d}{dt}\int_\Omega (w-k)_+^2 dx &\geq 2\int_\Omega \left(\nu|D_xw|^2 + \lambda(w-k)\int_0^1 D_pH(t,x,sD_xu,\mu(t)) ds \cdot D_xw\right)\chi_{w \geq k}dx \\
        &\geq -C_\nu\lambda^2\|D_pH\|_\infty^2\int_\Omega (w-k)_+^2 dx.        
    \end{aligned}
    $$
    By Gronwall's inequality,
    $$\int_\Omega (w(t,x)-k)_+^2 dx \leq e^{C_\nu\lambda^2\|D_pH\|_\infty^2(T-t)}\int_\Omega (w(T,x)-k)_+^2 dx.$$
    Choosing $k = \underset{x}{\max} \, w(T,x)_+ = \underset{x}{\max} \, u(T,x)_+$ gives $u(t,x) \leq \underset{x}{\max} \, u(T,x)_+ + \lambda\int_t^T \|H(t,x,0,\mu(s))\|_\infty ds$. By an analogous argument, $u(t,x) \geq -\underset{x}{\max} \, u(T,x)_- - \lambda\int_t^T \|H(t,x,0,\mu(s))\|_\infty ds$.
\end{proof}

\begin{lemma} \label{Lem: Bounding u in terms of Du}
    Assume A\ref{A: Hamiltonian}-\ref{A: Bound for DpH}, A\ref{A: Bounding H} and A\ref{A: DpH is Lipschitz in alpha} hold. Then for any $\theta \in (0,1)$, we have
    $$\|u\|_\infty \leq \|u(T,\cdot)\|_\infty + C_0\lambda T + \frac{\lambda_2\lambda^{q'+1}C_0^{q'}}{(1-\lambda\lambda_0)^{q'}}\left(\theta^{1-q'}T + (1-\theta)^{1-q'}\int_0^T\int_\Omega |D_xu|^q dm(t,x)dt\right).$$
    As a corollary,
    $$\|u\|_\infty \leq \|u(T,\cdot)\|_\infty + C_0\lambda T + \frac{\lambda_2\lambda^{q'+1}C_0^{q'}}{(1-\lambda\lambda_0)^{q'}}\left(\theta^{1-q'} + (1-\theta)^{1-q'}\|D_xu\|_\infty^q\right)T.$$
\end{lemma}
    
\begin{proof}
    By Proposition \ref{Prop: Maximum Principle}, we have
    $$\|u\|_\infty \leq \|u(T,\cdot)\|_\infty + \lambda\int_0^T\|H(t,\cdot,0,\mu(t))\|_\infty dt.$$
    By Lemma \ref{Lem: Existence & Bound for mu 1}, we get
    $$
    \begin{aligned}
        |H(t,x,0,\mu)| &\leq C_0 + \lambda_2\Lambda_{q_0}(\mu)^{q'} \\
        &\leq C_0 + \frac{\lambda_2\lambda^{q'}C_0^{q'}}{(1-\lambda\lambda_0)^{q'}}\left(1 + \||D_xu|^{q-1}\|_{L^{q_0}(m)}\right)^{q'} \\
        &\leq C_0 + \frac{\lambda_2\lambda^{q'}C_0^{q'}}{(1-\lambda\lambda_0)^{q'}}\left(\theta^{1-q'} + (1-\theta)^{1-q'}\||D_xu|^{q-1}\|_{L^{q_0}(m)}^{q'}\right).
    \end{aligned}
    $$
    Using $q_0 \leq q'$ and recalling $m(t,\cdot)$ is a probability density, we have $\||D_xu|^{q-1}\|_{L^{q_0}(m)}^{q'} \leq \|D_xu\|_{L^q(m)}^q \leq \|D_x u\|_{\infty}^q$, which concludes the proof.
\end{proof}
    
\begin{lemma} \label{Lem: Lq estimate for Du}
    Assume A\ref{A: Hamiltonian}-\ref{A: Convexity assumption for H} hold. Then
    \begin{multline*}
        \int_0^T\int_\Omega |D_xu|^q dm(t,x)dt
        \\
        \leq \left(1-\frac{\lambda_1\lambda^{q'}C_0^{q'}}{(1-\lambda\lambda_0)^{q'}(1-\theta)^{q'-1}}\right)^{-1}\left(C_0^2(1 + T) + C_0\|u\|_\infty + \frac{\lambda_1\lambda^{q'}C_0^{q'}}{(1-\lambda\lambda_0)^{q'}\theta^{q'-1}}T\right)
    \end{multline*}
    for all $\theta \in (0,1)$ with $\frac{\lambda_1\lambda^{q'}C_0^{q'}}{(1-\lambda\lambda_0)^{q'}}(1-\theta)^{1-q'} < 1$.
\end{lemma}

\begin{proof}
    Multiply the HJ equation in \eqref{Eq: Parametrized MFGC 1} by $m$ and integrate by parts to get the standard ``energy identity" for mean field games:
    $$\int_\Omega ( g(x,m(T,x))m(T,x) - u(0,x)m_0(x)) dx + \int_0^T\int_\Omega (D_pH \cdot D_xu - H) dm(t,x)dt = 0.$$
    By A\ref{A: Convexity assumption for H},
    $$
    \begin{aligned}
        &\int_0^T\int_\Omega |D_xu|^q dm(t,x)dt \\
        &\leq C_0^2T + C_0\int_\Omega (u(0,x)m_0(x) - g(x)m(T,x))dx + \lambda_1\int_0^T\Lambda_{q_0}(\mu)^{q'} dt \\
        &\leq C_0^2(1 + T) + C_0\|u\|_\infty + \frac{\lambda_1\lambda^{q'}C_0^{q'}}{(1-\lambda\lambda_0)^{q'}}\left(\theta^{1-q'}T + (1-\theta)^{1-q'}\int_0^T\int_\Omega |D_xu|^q dxdt\right).
    \end{aligned}
    $$
    Thus for $\theta \in (0,1)$ with $\frac{\lambda_1\lambda^{q'}C_0^{q'}}{(1-\lambda\lambda_0)^{q'}(1-\theta)^{q'-1}} < 1$, we have
    $$\int_0^T\int_\Omega |D_xu|^q dm(t,x)dt \leq \left(1-\frac{\lambda_1\lambda^{q'}C_0^{q'}}{(1-\lambda\lambda_0)^{q'}(1-\theta)^{q'-1}}\right)^{-1}\left(C_0^2(1 + T) + C_0\|u\|_\infty + \frac{\lambda_1\lambda^{q'}C_0^{q'}\theta^{1-q'}}{(1-\lambda\lambda_0)^{q'}}T\right).$$
\end{proof}

\begin{proposition}
    Assume A\ref{A: Hamiltonian}-\ref{A: Bounding H} and A\ref{A: DpH is Lipschitz in alpha} hold. Then $\|u\|_\infty \leq C$, where $C$ depends only on $\|u(T,\cdot)\|_\infty$ and the constants in the assumptions.
    \label{Prop: Bound for u}
\end{proposition}

\begin{proof}
    Choose $\theta \in (0,1)$ such that $\lambda_1\lambda^{q'} + C_0\lambda_2\lambda^{q'+1} < \frac{(1-\theta)^{q'-1}(1-\lambda\lambda_0)^{q'}}{C_0^{q'}}$. By lemmas \ref{Lem: Bounding u in terms of Du} and \ref{Lem: Lq estimate for Du}, we have
    $$
    \begin{aligned}
        \|u\|_\infty &\leq \|u(T,\cdot)\|_\infty + C_0\lambda T + \frac{\lambda_2\lambda^{q'+1}C_0^{q'}}{(1-\lambda\lambda_0)^{q'}}\left(\theta^{1-q'}T + (1-\theta)^{1-q'}\int_0^T\int_\Omega |D_xu|^q dm(t,x)dt\right) \\
        &\leq \frac{\lambda_2\lambda^{q'+1}C_0^{q'+1}}{(1-\lambda\lambda_0)^{q'}(1-\theta)^{q'-1}}\left(1-\frac{\lambda_1\lambda^{q'}C_0^{q'}}{(1-\lambda\lambda_0)^{q'}(1-\theta)^{q'-1}}\right)^{-1}\|u\|_\infty + C_\theta \\
        &= C_0\lambda_2\lambda^{q'+1}\left(\frac{(1-\lambda\lambda_0)^{q'}(1-\theta)^{q'-1}}{C_0^{q'}}-\lambda_1\lambda^{q'}\right)^{-1}\|u\|_\infty + C_\theta
    \end{aligned}
    $$
    Since $C_0\lambda_2\lambda^{q'+1}\left(\frac{(1-\lambda\lambda_0)^{q'}(1-\theta)^{q'-1}}{C_0^{q'}}-\lambda_1\lambda^{q'}\right)^{-1} < 1$, this concludes the proof.
\end{proof}

\subsection{Gradient Estimate} \label{sec: Gradient Bounds}

In this section we prove a priori bounds on $|D_x u|$ for solutions of System \eqref{Eq: Parametrized MFGC 1}.
We use what is commonly known as a Bernstein method.
For mean field games of controls, we follow an outline similar to that which is found in \cite{kobeissi2022classical}, where the boundary conditions play no role.
Here, we will need to adapt the argument to the case of Neumann or Dirichlet type boundary conditions.
We will start with the Neumann case.

\begin{theorem} \label{Thm: Gradient Estimate (N)}
    Assume A\ref{A: Hamiltonian}-\ref{A: DpH is Lipschitz in alpha} hold. Let $(u,m,\mu)$ be a classical solution of \eqref{Eq: Parametrized MFGC 1}-\eqref{Eq: Parametrized Neumann 1} with $u \in C^{1,3}([0,T] \times \Omega)$.
    Then
    $$|D_xu| \leq C\lambda^\frac{1}{2}.$$
\end{theorem}

\begin{proof}
    By vector calculus,
    \begin{multline*}
        -\frac{1}{2}\partial_t|D_xu|^2 - \frac{\nu}{2}\Delta |D_xu|^2 + \nu|D_{xx}^2u|^2 + \frac{\lambda}{2}D_x|D_xu|^2 \cdot D_pH(t,x,D_xu,\mu)  \\
        = -\lambda D_xH(x,D_xu,\mu) \cdot D_xu.
    \end{multline*}
    Define the following functions:
    $$\varphi(v) = \exp\left(\exp\left(-a-b\|u\|_\infty^{-1} v\right)\right), \text{ for } |v| \leq \|u\|_\infty,$$
    $$w(t,x) = \varphi(u(T-t,x))|D_xu(T-t,x)|^2$$
    where $a > 1$ and $b > 0$ are constants that will be defined below. Note that
    $$\varphi'(v) = -b\|u\|_\infty^{-1}e^{-a-b\|u\|_\infty^{-1}v}\varphi(v),$$
    $$\varphi''(v) = b^2\|u\|_\infty^{-2}e^{-a-b\|u\|_\infty^{-1}v}\left(1+e^{-a-b\|u\|_\infty^{-1}v}\right)\varphi(v),$$
    and hence
    $$1 \leq \varphi(v) \leq e^{e^{-a+b}}, \hspace{2cm} b\|u\|_\infty^{-1}e^{-a-b} \leq \frac{|\varphi'(v)|}{\varphi(v)} \leq b\|u\|_\infty^{-1}e^{-a+b}.$$
    Furthermore,
    $$
    \begin{aligned}
        &\partial_tw - \nu\Delta w + \lambda D_xw \cdot D_pH(t,x,D_xu,\mu) + 2\nu\frac{\varphi'}{\varphi}D_xw \cdot D_xu + 2\nu\varphi|D_{xx}^2 u|^2 \\
        &= \frac{\varphi'}{\varphi}w\left[-\partial_tu - \nu\Delta u + \lambda D_xu \cdot D_pH(t,x,D_xu,\mu)\right] - \nu\frac{\varphi''\varphi-2(\varphi')^2}{\varphi^3}w^2 - 2\lambda\varphi D_xu \cdot D_xH(x,D_xu,\mu) \\
        &= \lambda \frac{\varphi'}{\varphi}w\left[D_xu \cdot D_pH(t,x,D_xu,\mu) - H(t,x,D_xu,\mu)\right] - \nu\frac{\varphi''\varphi-2(\varphi')^2}{\varphi^3}w^2 - 2\lambda\varphi D_xu \cdot D_xH(x,D_xu,\mu).
    \end{aligned}
    $$
    We now bound the right-hand side from above. To begin, notice that $\varphi''\varphi - 2(\varphi')^2 \geq 0$ provided $a \geq b$. By A\ref{A: Bound for DxH},
    $$-2\varphi D_xu \cdot D_xH(x,D_xu,\mu) \leq 2C_0\varphi|D_xu|\left(1 + |D_xu|^q + \Lambda_{q_0}(\mu)^{q'}\right)$$
    Since $\varphi' < 0$, A\ref{A: Convexity assumption for H} gives
    $$
    \begin{aligned}
        \frac{\varphi'}{\varphi}w\left[D_xu \cdot D_pH(t,x,D_xu,\mu) - H(t,x,D_xu,\mu)\right] \leq -C_0^{-1}\frac{|\varphi'|}{\varphi^{1+\frac{q}{2}}}w^{1+\frac{q}{2}} + C_0\frac{|\varphi'|}{\varphi}w + C_0^{-1}\lambda_1\frac{|\varphi'|}{\varphi}w\Lambda_{q_0}(\mu)^{q'}.
    \end{aligned}
    $$
    By Lemma \ref{Lem: Existence & Bound for mu 1}, we get
    $$\frac{|\varphi'|}{\varphi}\Lambda_{q_0}(\mu)^{q'} \leq b\|u\|_\infty^{-1}e^{-a+b}\frac{C_0^{q'}}{(1-\lambda_0)^{q'}}\left(\theta^{1-q'} + (1-\theta)^{1-q'}\|w\|_\infty^{\frac{q}{2}}\right)$$
    for $\theta \in (0,1)$ defined below. Also note that
    $$\frac{|\varphi'|}{\varphi^{1+\frac{q}{2}}} \geq b\|u\|_\infty^{-1}e^{-a-b}e^{-\frac{q}{2}e^{-a+b}}.$$
    Combining these inequalities, we get
    \begin{equation}
        \begin{aligned}
            \partial_tw - \nu\Delta w + \lambda D_xw \cdot D_pH(t,x,D_xu,\mu) - 2\nu\frac{\varphi'}{\varphi}D_xw \cdot D_xu + 2\nu\varphi|D_{xx}^2 u|^2 \\
            \leq \lambda\left[-\Tilde{C}w^{1+\frac{q}{2}} + b\|u\|_\infty^{-1}e^{-a+b}\frac{\lambda_1C_0^{q'-1}}{(1-\theta)^{1-q'}(1-\lambda_0)^{q'}}\|w\|_\infty^{1+\frac{q}{2}} + C_{a,b,\theta}(1 + \|u\|_\infty^{-1})(1+\|w\|_\infty^{1+\frac{q}{2}})\right]
        \end{aligned}
        \label{Eq: w inequality 1}
    \end{equation}
    where $\Tilde{C} = C_0^{-1}b\|u\|_\infty^{-1}e^{-a-b}e^{-\frac{q}{2}e^{-a+b}}$. Since $\|w(0,\cdot)\|_\infty \leq eC_0^2$, we get that the constant function
    $$\Tilde{v} = \lambda^\frac{2}{2+q}\left[\frac{\lambda_1C_0^{q'}e^{2b}e^{\frac{q}{2}e^{-a+b}}}{(1-\theta)^{q'-1}(1-\lambda_0)^{q'}}\|w\|_\infty^{1+\frac{q}{2}} + C(1+\|u\|_\infty)\left(1+\|w\|_\infty^{\frac{1+q}{2}}\right)\right]^{\frac{2}{2+q}}$$
    is a super-solution to \eqref{Eq: w inequality 1} with $\|w(0,\cdot)\|_\infty \leq \Tilde{v}$ (where $C = C + (eC_0^2)^{1+\frac{q}{2}}$). Note that $v = w - \Tilde{v}$ satisfies
    $$\partial_tv - \nu \Delta v + D_xv \cdot \left(\lambda D_pH - 2\nu\frac{\varphi'}{\varphi}D_xu\right) \leq \Tilde{C}\left(\Tilde{v}^{1+\frac{q}{2}} - w^{1+\frac{q}{2}}\right).$$
    
    Recall that we are assuming in the case of Neumann boundary conditions that $\Omega$ is convex. Thus, for all $z \in C^2(\overline{\Omega})$ with $\frac{\partial z}{\partial\n} = 0$ on $\B$, we have $\frac{\partial}{\partial\n}|D_xz|^2 \leq 0$ on $\B$ (See \cites{lions1980resolution,lions1985quelques}). Hence, we get $\frac{\partial}{\partial\n}w \leq 0$. Since
    $$\frac{d}{dt}\int_\Omega v_+^2 dx + 2\nu\int_\Omega |D_xv|^2\chi_{w \geq \Tilde{v}} dx \leq 2\nu\int_\B v_+D_xv \cdot \n d\sigma(x) + 2\int_\Omega v_+D_xv \cdot \left(\lambda D_pH - 2\nu\frac{\varphi'}{\varphi}D_xu\right) dx,$$
    this implies
    $$\frac{d}{dt}\int_\Omega v_+^2 dx \leq C_{\nu,\lambda}\left(\lambda ^2\|D_pH\|_\infty^2 + 4\nu^2\left\|\frac{\varphi'}{\varphi}\right\|_\infty^2\|D_xu\|_\infty^2\right)\int_\Omega v_+^2 dx.$$
    By Gronwall's inequality, it follows that $w \leq \Tilde{v}$. This implies
    $$\|w\|_\infty^{1+\frac{q}{2}} \leq \lambda\left[\frac{\lambda_1C_0^{q'}e^{2b}e^{\frac{q}{2}e^{-a+b}}}{(1-\theta)^{q'-1}(1-\lambda_0)^{q'}}\|w\|_\infty^{1+\frac{q}{2}} + C(1+\|u\|_\infty)\left(1+\|w\|_\infty^{\frac{1+q}{2}}\right)\right]$$
    By A\ref{A: Convexity assumption for H}, we can choose $a,b,\theta$ such that $\frac{\lambda_1C_0^{q'}e^{2b}e^{\frac{q}{2}e^{-a+b}}}{(1-\theta)^{q'-1}(1-\lambda_0)^{q'}} < 1$ and so
    $$\|w\|_\infty^{1+\frac{q}{2}} \leq \lambda C'(1+\|u\|_\infty)\left(1+\|w\|_\infty^{\frac{1+q}{2}}\right).$$
    If $\|w\|_\infty^\frac{1}{2} \leq 2\lambda^\frac{1}{2}C'(1+\|u\|_\infty)$, we are done. So suppose $\|w\|_\infty^\frac{1}{2} > 2\lambda^\frac{1}{2}C'(1+\|u\|_\infty)$. Then we have $\|w\|_\infty^{1+\frac{q}{2}} \leq \frac{1}{2}\lambda^\frac{1}{2}\|w\|_\infty^{\frac{1}{2}}\left(1+\|w\|_\infty^{\frac{1+q}{2}}\right)$, which implies $\|w\|_\infty \leq \lambda^\frac{1}{2}$. Since $\varphi \geq 1$, this completes the proof.
\end{proof}

To handle the Dirichlet case, we first need to prove an a priori bound on $\abs{D_x u}$ on the boundary.
Arguments of a similar spirit can be found in \cite[Chapter V]{ladyzhenskaia1968linear}.
Here it is necessary to have $q \leq 2$ (see Assumption \ref{A: Bound for DpH}).
Indeed, for superquadratic Hamiltonians, the phenomenon of ``gradient blow-up" is well-known for Hamilton-Jacobi equations with Dirichlet boundary conditions.
See e.g.~\cite{attouchi2020gradient} and references therein.
\begin{lemma} \label{Lem: Bounds for Du on Boundary (D)}
    Assume A\ref{A: Hamiltonian}-\ref{A: DpH is Lipschitz in alpha} hold. Let $(u,m,\mu)$ be a classical solution of \eqref{Eq: Parametrized MFGC 1}-\eqref{Eq: Parametrized Dirchlet 1}.
    Then
    $$|D_xu| \leq C$$
    on $\Sigma$.
\end{lemma}

\begin{proof}
    Let $S$ be a section of $\B$ such that, under a $C^2$ smooth change of coordinates $y = y(x)$, the image of $S$ is contained in $\R^{n-1} \times \{0\}$ and the image of $\overline{\Omega}$ is contained in $\R^n_+ = \{y \in \R^n : y_n \geq 0\}$. We will denote by $\Tilde{S}$ and $\Tilde{\Omega}$ the images of $S$ and $\Omega$, respectively. Define $\Tilde{u}(t,y(x)) = u(t,x)$. Note that
    $$\partial_iu(t,x) = C_{ik}\partial_k\Tilde{u}(t,y), \hspace{2cm} \partial_{ij}u(t,x) = A_{ijkl}\partial_{lk}\Tilde{u}(t,y) + B_{ijk}\partial_k\Tilde{u}(t,y)$$
    where we follow the convention of summing over repeated indices, and we define
    $$A_{ijkl} = \partial_iy_k\partial_jy_l, \hspace{2cm} B_{ijk} = \partial_{ij}y_k, \hspace{2cm} C_{ij} = \partial_iy_j.$$
    Let $A_{kl} = \sum_{i=1}^n A_{iikl}$, $B_k = \sum_{i=1}^n B_{iik}$, and $\Tilde{H}(t,y,D_y\Tilde{u}(t,y)\mu) = \lambda H(t,x,Du(t,x),\mu) - \nu B_k\partial_k\Tilde{u}(t,y)$. Note that all coefficients are bounded, there exists $a_0 > 0$ with
    $$A_{kl}\xi_k\xi_l \geq a_0|\xi|^2, \hspace{1cm} \forall\xi \in \R^n,$$
    and
    $$\Tilde{H}(t,y,D_y\Tilde{u}(t,y),\mu) \leq \Tilde{M}(1 + |D_y\Tilde{u}(t,y)| + |D_y\Tilde{u}(t,y)|^q) \leq M(1 + |D_y\Tilde{u}(t,y)|^2)$$
    where $M$ depends on the change of coordinates. (Here we have used $q \leq 2$.) Now note that $\Tilde{u}$ satisfies
    $$\begin{cases}
        -\partial_t\Tilde{u} - \nu A_{kl}\partial_{kl}\Tilde{u} + \Tilde{H}(t,y,D_y\Tilde{u},\mu) = 0, &(t,x) \in (0,T) \times \Tilde{\Omega} \\
        \Tilde{u}(t,x) = 0, &(t,x) \in [0,T] \times \partial\Tilde{\Omega}
    \end{cases}$$
    Define $v = e^{\gamma\Tilde{u}} - 1$ with $\gamma > 0$ to be chosen below. Then
    $$-\partial_tv - \nu A_{kl}\partial_{kl}v + \nu\gamma^2e^{\gamma\Tilde{u}}A_{kl}\partial_k\Tilde{u}\partial_l\Tilde{u} + \gamma e^{\gamma\Tilde{u}}\Tilde{H}(t,y,D_y\Tilde{u},\mu) = 0$$
    and so
    $$
    \begin{aligned}
        -\partial_tv - \nu A_{kl}\partial_{kl}v + a_0\nu\gamma^2e^{\gamma\Tilde{u}}|D_x\Tilde{u}|^2 &\leq M\gamma(1 + |D_y\Tilde{u}|^2)e^{\gamma\Tilde{u}}
    \end{aligned}
    $$
    Choosing $\gamma = \frac{M}{a_0\nu}$ gives
    $$-\partial_tv - \nu A_{kl}\partial_{kl}v \leq \frac{M^2}{a_0\nu}e^{\gamma\Tilde{u}} \leq C,$$
    using the bound on $\tilde u$ coming from the maximum principle.
    
    Define $w(t,y) = v(t,y) + \eta e^{-y_n}$, where $\eta$ is a constant to be chosen below. Note that
    $$-\partial_tw - \nu A_{kl}\partial_{kl}w \leq C - \nu\eta A_{nn}e^{-d}$$
    where $d$ is the diameter of $\Tilde{\Omega}$. For $\eta$ large enough, we get
    $$-\partial_tw - \nu A_{kl}\partial_{kl}w < 0.$$
    By the maximum principle,
    $$w(t,y) \leq \underset{([0,T] \times \partial\Tilde{\Omega}) \cup (\{T\} \times \Tilde{\Omega)}}{\max} w.$$
    Since $v(t,y) = 0$ for all $y \in \partial\Tilde{\Omega}$ and since $\Tilde{\Omega} \subseteq \R^n_+$, it follows that
    $$\underset{[0,T] \times \partial\Tilde{\Omega}}{\max} w = \eta.$$
    Moreover, for $\eta$ large enough, we have
    $$\partial_nw(T,y) = \gamma e^{\gamma\Tilde{u}(T,y)}\partial_n\Tilde{u}(T,y) - \eta e^{-y_n} \leq \gamma e^{\gamma\|g\|_\infty}\|C_{ik}\|_\infty\|D_xg\|_\infty - \eta e^{-d} \leq 0$$
    and hence
    $$\max_{\Tilde{\Omega}} w(T,y) = \max_{\partial\Tilde{\Omega}} = \eta.$$
    Since $w$ attains its maximum value of $\eta$ on $[0,T] \times \Tilde{S}$ provided $\eta$ is large enough, we get that for $y \in \Tilde{S}$, $\partial_nw \leq 0$ and hence $\partial_n\Tilde{u} \leq \frac{1}{\gamma}\eta e^{\gamma\|u\|_\infty}$. Since $D_xu \cdot \n = C\partial_n\Tilde{u}$, this gives an upper bound on the normal derivative of $u$ on the boundary. A similar argument gives a lower bound on the normal derivative. Since $u = 0$ on $\Sigma$, the tangential derivatives are all $0$.
\end{proof}

\begin{theorem} \label{Thm: Gradient Estimate (D)}
    Assume A\ref{A: Hamiltonian}-\ref{A: DpH is Lipschitz in alpha} hold. Let $(u,m,\mu)$ be a classical solution of \eqref{Eq: Parametrized MFGC 1}-\eqref{Eq: Parametrized Dirchlet 1} with $u \in C^{1,3}([0,T] \times \Omega)$. Then
    $$|D_xu| \leq C.$$
\end{theorem}

\begin{proof}
    The argument is almost identical to the proof of theorem \ref{Thm: Gradient Estimate (N)}, taking into account Lemma \ref{Lem: Bounds for Du on Boundary (D)} to estimate the gradient on the boundary.
\end{proof}

\subsection{Bootstrapping}
\label{Sec: Bootstrapping 1}

Let $X \coloneqq C^{\beta/2,1+\beta}(Q) \times C^0(0,T;L^2(\Omega))$. With the aim of applying the Leray-Schauder fixed point theorem, we will define map the $\Psi: X \times [0,1] \rightarrow X$ as follows: Given $(\Tilde{u},\Tilde{m}) \in X$ and $\lambda \in [0,1]$, define
$$\overline{m}(t,x) =
\begin{cases}
    \frac{|\Tilde{m}(t,x)|}{\|\Tilde{m}(t,\cdot)\|_{L^1(\Omega)}}, &\text{if } \|\Tilde{m}(t,\cdot)\|_{L^1(\Omega)} > 1 \\
    |\Tilde{m}(t,x)|, &\text{if } 0 < \|\Tilde{m}(t,\cdot)\|_{L^1(\Omega)} \leq 1 \\
    0, &\text{if } \|\Tilde{m}(t,\cdot)\|_{L^1(\Omega)} = 0
\end{cases}$$
and let $(u,m) = \Psi(\Tilde{u},\Tilde{m},\lambda)$ be the classical solution to
\begin{equation}
    \begin{cases}
        -u_t - \nu\Delta u + \lambda H(t,x,D_x\Tilde{u},\overline{\mu}) = 0, &(t,x) \in Q \\
        m_t - \nu\Delta m - \lambda\nabla \cdot (mD_pH(t,x,D_x\Tilde{u},\Tilde{\mu})) = 0, &(t,x) \in Q \\
        \Tilde{\mu} = (I,-\lambda D_pH(t,\cdot,D_x\Tilde{u},\Tilde{\mu})) \# \overline{m}, &t \in [0,T] \\
        \overline{\mu} = (I,-\lambda D_pH(t,\cdot,D_x\Tilde{u},\overline{\mu})) \# m, &t \in [0,T] \\
        u(T,x) = \lambda g(x,m(T)), \hspace{1cm} m(0,x) = \lambda m_0(x), &x \in \overline{\Omega} \\
    \end{cases}
    \label{Eq: Definition of Psi}
\end{equation}
paired with either
\begin{equation}
    u = m = 0, \hspace{1cm} (t,x) \in \Sigma \tag{\ref{Eq: Definition of Psi}d}
    \label{Eq: Psi Dirchlet}
\end{equation}
or
\begin{equation}
    \frac{\partial u}{\partial\n} = \nu\frac{\partial m}{\partial\n} + \lambda mD_pH(t,x,D_x\Tilde{u},\Tilde{\mu}) \cdot \n = 0, \hspace{1cm} (t,x) \in \Sigma \tag{\ref{Eq: Definition of Psi}n}
    \label{Eq: Psi Neumann}
\end{equation}

\begin{remark}
    By A\ref{A: Bound for DpH}, A\ref{A: Bounding H}, A\ref{A: Bound for DxH}, we have uniform bounds for $H$ and each of its first-order derivatives (depending only on $\|D_x\Tilde{u}\|_\infty$ and the constants in the assumptions). By Lemma \ref{Lem: Existence & Bound for mu 1} and the arguments from \cite[Section 7.1.2]{evans2022partial}, we get that $\Tilde{\mu}$, $m$, and $\overline{\mu}$ are well-defined.
\end{remark}

In order to prove the continuity and compactness required to apply Leray-Schauder, we will need prove higher regularity and a priori estimates for $(u,m)$, possibly depending on $\|\Tilde{u}\|_{C^{\beta/2,1+\beta}(Q)}$ and $\|\Tilde{m}\|_{C^0(0,T;L^2(\Omega)}$.

\begin{lemma} \label{Lem: Boundedness of m}
    Assume A\ref{A: Hamiltonian}-\ref{A: DpH is Lipschitz in alpha} hold. Then there exists $C > 0$ (depending only on $\|D_x\Tilde{u}\|_\infty$, $\|m_0\|_\infty$, $\nu$, $n$, $T$, and $|\Omega|$) such that
    $$\|m\|_{L^\infty(Q)} \leq C.$$
\end{lemma}

\begin{proof}
    In the Dirichlet case, this follows directly from \cite[Theorem V.2.1]{ladyzhenskaia1968linear}. For the convenience of the reader, we sketch a proof that holds in both Dirichlet and Neumann problems. We use Moser iteration. First note that for $p \geq 2$, we have
    $$
    \begin{aligned}
        &\int_\Omega m^p dx + \nu p(p-1)\int_0^t\int_\Omega m^{p-2}|D_xm|^2 dxd\tau \\
        &= \lambda\int_\Omega m_0^p dx - \lambda p(p-1)\int_0^t\int_\Omega m^{p-1}D_pH \cdot D_xm dxd\tau \\
        &\leq \int_\Omega m_0^p dx + \frac{\nu p(p-1)}{2}\int_0^t\int_\Omega m^{p-2}|D_xm|^2 dxd\tau + \frac{p(p-1)}{2\nu}\|D_pH\|_\infty^2\int_0^t\int_\Omega m^p dxd\tau.
    \end{aligned}
    $$
    The inequality is established formally by multiplying the FP equation by $m^{p-1}$, which is a valid test function when $m$ is bounded; to establish the inequality for all weak solutions, one can multiply by $\phi(m)$ where $\phi(s)$ is a smooth bounded function chosen to approximate $s^{p-1}$ for $s \geq 0$ (we omit the details).
    Now, by Sobolev's inequality,
    $$\|m(t)\|_{L^\frac{pn}{n-2}(\Omega)} = \|m(t)^\frac{p}{2}\|_{L^{2^*}(\Omega)}^\frac{2}{p} \leq C_n^\frac{2}{p}(\|m(t)\|_{L^p(\Omega)}^p + \|D_x(m(t)^\frac{p}{2})\|_{L^2(\Omega)}^2)^\frac{1}{p}$$
    and so
    $$
    \begin{aligned}
        &\|m\|_{L^p\left(0,T;L^\frac{pn}{n-2}(\Omega)\right)} \\
        &\leq C_n^\frac{1}{p}\left(\int_0^T (\|m(t)\|_{L^p(\Omega)}^p + \|D_x(m(t)^\frac{p}{2})\|_{L^2(\Omega)}^2) dt\right)^\frac{1}{p} \\
        &\leq C_n^\frac{1}{p}\left(T + \frac{p}{2\nu(p-1)}\right)^\frac{1}{p}\left(\|m_0\|_{L^p(\Omega)}^p + \frac{p(p-1)}{2\nu}\|D_pH\|_\infty^2\|m\|_{L^p(Q)}^p\right)^\frac{1}{p}.
    \end{aligned}
    $$
    By interpolation,
    $$
    \begin{aligned}
        \int_\Omega m^{(1 + \frac{2}{n})p} dx &\leq \left(\int_\Omega m^\frac{pn}{n-2} dx\right)^\frac{n-2}{n}\left(\int_\Omega m^p dx\right)^\frac{2}{n} \\
        &\leq \left(\int_\Omega m^\frac{pn}{n-2} dx\right)^\frac{n-2}{n}\left(\int_\Omega m_0^p dx + \frac{p(p-1)}{2\nu}\|D_pH\|_\infty^2\int_0^T\int_\Omega m^p dxdt\right)^\frac{2}{n}.
    \end{aligned}
    $$
    Thus,
    $$
    \begin{aligned}
        &\|m\|_{L^{(1 + \frac{2}{n})p}(Q)} \\
        &\leq \|m\|_{L^p\left(0,T;L^\frac{pn}{n-2}(\Omega)\right)}^\frac{1}{1 + \frac{2}{n}}\left(\|m_0\|_{L^p(\Omega)}^p + \frac{p(p-1)}{2\nu}\|D_pH\|_\infty^2\|m\|_{L^p(Q)}^p\right)^\frac{\frac{2}{n}}{(1 + \frac{2}{n})p} \\
        &\leq C_n^\frac{1}{(1 + \frac{2}{n})p}\left(T + \frac{p}{2\nu(p-1)}\right)^\frac{1}{(1 + \frac{2}{n})p}\left(\|m_0\|_{L^p(\Omega)}^p + \frac{p(p-1)}{2\nu}\|D_pH\|_\infty^2\|m\|_{L^p(Q)}^p\right)^\frac{1}{p} \\
        &\leq C_n^\frac{1}{(1 + \frac{2}{n})p}\left(T + \frac{p}{2\nu(p-1)}\right)^\frac{1}{(1 + \frac{2}{n})p}\left(|\Omega| + \frac{p(p-1)}{2\nu}\|D_pH\|_\infty^2\right)^\frac{1}{p}\max\{\|m_0\|_{L^\infty(\Omega)},\|m\|_{L^p(Q)}\}.
    \end{aligned}
    $$
    Let $\gamma = 1 + \frac{2}{n}$ and define a sequence $(p_k)_{k=0}^\infty$ by $p_k = 2\gamma^k$. Now define the sequence $(a_k)_{k=0}^\infty$ by $a_k = \max\{\|m_0\|_{L^\infty(\Omega)},\|m\|_{L^{p_k}(Q)}\}$. By induction,
    $$a_{k+1} \leq a_0\prod_{j=0}^k C_n^\frac{1}{p_{j+1}}\left(T + \frac{p_j}{2\nu(p_j-1)}\right)^\frac{1}{p_{j+1}}\left(|\Omega| + \frac{p_j(p_j-1)}{2\nu}\|D_pH\|_\infty^2\right)^\frac{1}{p_j}$$
    As the series
    $$
    \begin{aligned}
        &\sum_{j=0}^\infty \left[\frac{1}{p_{j+1}}\ln{C_n} + \frac{1}{p_{j+1}}\ln{\left(T + \frac{p_j}{2\nu(p_j-1)}\right)} + \frac{1}{p_j}\ln{\left(|\Omega| + \frac{p_j(p_j-1)}{2\nu}\|D_pH\|_\infty^2\right)}\right] \\
        &= \sum_{j=0}^\infty \left[\frac{1}{2\gamma^{j+1}}\ln{C_n} + \frac{1}{2\gamma^{j+1}}\ln{\left(T + \frac{2\gamma^j}{2\nu(2\gamma^j-1)}\right)} + \frac{1}{2\gamma^j}\ln{\left(|\Omega| + \frac{2\gamma^j(2\gamma^j-1)}{2\nu}\|D_pH\|_\infty^2\right)}\right]
    \end{aligned}
    $$
    converges, it follows that
    $$\prod_{j=0}^\infty C_n^\frac{1}{p_{j+1}}\left(T + \frac{p_j}{2\nu(p_j-1)}\right)^\frac{1}{p_{j+1}}\left(|\Omega| + \frac{p_j(p_j-1)}{2\nu}\|D_pH\|_\infty^2\right)^\frac{1}{p_j}$$
    also converges, thus completing the proof.
\end{proof}

\begin{lemma} \label{Lem: Holder Continuity of m}
    Assume A\ref{A: Hamiltonian}-\ref{A: DpH is Lipschitz in alpha} hold. Then there exist $\alpha \in (0,\beta)$ and $C > 0$ so that $m \in C^{\alpha/2,\alpha}(Q)$ with
    $$\|m\|_{C^{\alpha/2,\alpha}(Q)} \leq C,$$
    where $C$ depends on $(\Tilde{u},\Tilde{m})$.
\end{lemma}

\begin{proof}
    In the Dirichlet case, this follows from \cite[Theorem V.1.1]{ladyzhenskaia1968linear} or \cite[Theorem II.1.2]{dibenedetto2012degenerate}. In the Neumann case, this follows by arguments similar to those used to prove \cite[Theorem II.1.3]{dibenedetto2012degenerate}.
\end{proof}

\begin{lemma}
    Assume A\ref{A: Hamiltonian}-\ref{A: DpH is Lipschitz in m} hold. Then $h(t,x) \coloneqq D_pH(t,x,D_x\Tilde{u}(t,x),\overline{\mu}(t)) \in C^{\alpha\beta/2,\alpha\beta}(Q)$ with$$\|h\|_{C^{\alpha\beta/2,\alpha\beta}(Q)} \leq C(\|D_x\Tilde{u}\|_{ C^{\alpha/2,\alpha}(Q)}^\beta + \|m\|_{C^{\alpha/2,\alpha}(Q)}^\beta + 1).$$
\end{lemma}

\begin{proof}
    Note that for $x,y \in \overline{\Omega}$ and $t \in [0,T]$,
    $$
    \begin{aligned}
        |h(t,x) - h(t,y)| \leq& |D_pH(t,x,D_x\Tilde{u}(t,x),\overline{\mu}(t)) - D_pH(t,y,D_x\Tilde{u}(t,x),\overline{\mu}(t))| \\
        &+ |D_pH(t,y,D_x\Tilde{u}(t,x),\overline{\mu}(t)) - D_pH(t,y,D_x\Tilde{u}(t,y),\overline{\mu}(t))| \\
        \leq& C(1 + \|D_x\Tilde{u}\|_{C^{\alpha/2,\alpha}(Q)}^\beta)|x-y|^{\alpha\beta}.
    \end{aligned}
    $$
    Furthermore, for all $t,s \in [0,T]$,
    $$
    \begin{aligned}
        \|h(t,\cdot) - h(s,\cdot)\|_\infty \leq& \|D_pH(t,\cdot,D_x\Tilde{u}(t,\cdot),\overline{\mu}(t)) - D_pH(t,\cdot,D_x\Tilde{u}(t,\cdot),(I,h(t,\cdot)) \# m(s))\|_\infty \\
        &+ \|D_pH(t,\cdot,D_x\Tilde{u}(t,\cdot),(I,h(t,\cdot)) \# m(s)) - D_pH(t,\cdot,D_x\Tilde{u}(t,\cdot),\overline{\mu}(s))\|_\infty \\
        &+ \|D_pH(t,\cdot,D_x\Tilde{u}(t,\cdot),\overline{\mu}(s)) - D_pH(s,\cdot,D_x\Tilde{u}(s,\cdot),\overline{\mu}(s))\|_\infty \\
        \leq& Cd^*(m(t),m(s))^\beta + L_1\|h(t,\cdot) - h(s,\cdot)\|_\infty \\
        &+ C(1 + \|D_x\Tilde{u}\|_{C^{\alpha/2,\alpha}(Q)}^\beta)|t-s|^{\alpha\beta/2}
    \end{aligned}
    $$
    and
    $$
    \begin{aligned}
        d^*(m(t),m(s)) &\leq \max\{1,diam(\Omega)\}\int_\Omega |m(t) - m(s)| dx \\
        &\leq \max\{1,diam(\Omega)\}|\Omega|\|m\|_{C^{\alpha/2,\alpha}(Q)}|t-s|^\frac{\alpha}{2},
    \end{aligned}
    $$
    As $L_1 \in (0,1)$, this completes the proof.
\end{proof}

\begin{remark}
    This gives us estimates for $\overline{\mu}$ in $C^\frac{\alpha\beta}{2}(0,T;\mathfrak{M}_{\infty,R}(\overline{\Omega} \times \R^n))$ (resp.~in $C^\frac{\alpha\beta}{2}(0,T;\mathcal{P}_{\infty,R}(\overline{\Omega} \times \R^n))$). By A\ref{A: Regularity of H} and A\ref{A: Regularity of DxH}, this gives us estimates for $H$ in $C^{\alpha\beta^2/2,\alpha\beta^2}(Q)$ and $D_xH$ in $C^{\alpha\beta^2/2,\alpha\beta^2}(Q;\R^n)$. By \cite[Theorem IV.5.2]{ladyzhenskaia1968linear} (resp., \cite[Theorem IV.5.3]{ladyzhenskaia1968linear}), this shows that $\Psi$ is well-defined. Furthermore, we get bounds for $u$ in $W^{1,2}_p(Q)$ for arbitrarily large $p$ (see \cite[Section IV.9]{ladyzhenskaia1968linear}).
    By \cite[Lemma II.3.3]{ladyzhenskaia1968linear}, we get estimates for $u$ and $D_xu$ in $C^{\beta/2,\beta}(Q)$ and $C^{\beta/2,\beta}(Q;\R^n)$, respectively.
\end{remark}

\begin{lemma} \label{Lem: u C1+a}
    Assume A\ref{A: Hamiltonian}-\ref{A: Regularity of H} hold. Then $u \in C^{1+\alpha\beta^2/2,2+\alpha\beta^2}(Q)$ with
    $$\|u\|_{C^{1+\alpha\beta^2/2,2+\alpha\beta^2}(Q)} \leq C\lambda(\|H\|_{C^{\alpha\beta^2/2,\alpha\beta^2}(Q)} + \|g\|_{C^{2+\alpha\beta^2}(\Omega)}).$$
\end{lemma}

\begin{proof}
    This follows from Theorem 5.2 (resp., Theorem 5.3) from chapter IV of \cite{ladyzhenskaia1968linear}.
\end{proof}

Finally, in order to use the results from Section \ref{sec: Gradient Bounds} to get a priori estimates for fixed-points, we will need to prove even higher of $u$.

\begin{lemma}
Assume A\ref{A: Hamiltonian}-\ref{A: Regularity of DxH} hold and fix $\psi \in C^{2+\alpha\beta^2}(\overline{\Omega})$ such that $\psi|_\B = 0$. 
Then for $j=1,\dots,n$, we get $\|\psi\partial_ju\|_{C^{1+\alpha\beta/2,2+\alpha\beta}(Q)} \leq C_\psi$. In particular, $u \in C^{1,3}([0,T] \times \Omega)$.
\end{lemma}

\begin{proof}
    Note that for $G \coloneqq \left(2\nu D_x\psi \cdot D_x\partial_ju + \nu\partial_ju\Delta\psi + \lambda\psi\left(\partial_jH + D_pH \cdot D_x\partial_j\Tilde{u}\right)\right)$, we have $G \in C^{\alpha\beta^2/2,\alpha\beta^2}(Q)$. By Theorem IV.5.2 from \cite{ladyzhenskaia1968linear}, there exists a solution $w \in C^{1+\alpha\beta^2/2,2+\alpha\beta^2}(Q)$ to
    $$
    \begin{cases}
        -\partial_tw - \nu\Delta w + G = 0, &x \in Q \\
        w(T,x) = \lambda\psi\partial_jg(x,m(T)), &x \in \overline{\Omega} \\
        w = 0, &x \in \Sigma
    \end{cases}
    $$
    and it satisfies
    $$\|w\|_{C^{1+\alpha\beta^2/2,2+\alpha\beta^2}(Q)} \leq C(\|G\|_{C^{\alpha\beta^2/2,\alpha\beta^2}(Q)} + \lambda\|\psi\partial_jg\|_{C^{2+\alpha\beta^2}(\Omega)}) < \infty.$$
    By uniqueness of solutions in the sense of distributions, it follows that $w = \psi\partial_ju$.
\end{proof}

\subsection{Existence}

With the results from Section \ref{Sec: Bootstrapping 1}, we are ready to prove our first existence result.

\begin{theorem} \label{Thm: Existence 1}
    Under assumptions A\ref{A: Hamiltonian}-\ref{A: Regularity of DxH}, there exists a classical solution to \eqref{Eq: MFGC}-\eqref{Eq: Dirchlet} (resp., \eqref{Eq: MFGC}-\eqref{Eq: Neumann}).
\end{theorem}

\begin{proof}
    Define $X$ and $\Psi: X \times [0,1] \rightarrow X$ as in section \ref{Sec: Bootstrapping 1}.

    \textit{$\Psi(\cdot,0)$ is constant:} 
    First, note that for every $(\Tilde{u},\Tilde{m}) \in X$, $\Psi(\Tilde{u},\Tilde{m},0) = (0,0)$.
    
    \textit{Bound for fixed-points:} 
    By the results from the previous sections, we get
    $$\|u\|_{C^{1+\alpha\beta^2/2,2+\alpha\beta^2}(Q)} + \|m\|_{C^0(0,T;L^2(\Omega))} \leq C$$
    for all $(u,m,\lambda) \in X \times [0,1]$ with $\Psi(u,m,\lambda) = (u,m)$.

    \textit{Continuity:} That $\Psi$ is continuous in $\lambda$ is clear. Now fix $\lambda \in [0,1]$. Take $(\Tilde{u}_k,\Tilde{m}_k) \rightarrow (\Tilde{u},\Tilde{m})$ in $X$. Define $(u_k,m_k) = \Psi(\Tilde{u}_k,\Tilde{m}_k,\lambda)$ and $(u,m) = \Psi(\Tilde{u},\Tilde{m},\lambda)$. By Remark \ref{rmk: uniform convergence of mu 1}, it follows that $\Tilde{\mu}^k = (I,-\lambda D_pH(t,\cdot,D_x\Tilde{u}_k,\Tilde{\mu}^k)) \# \overline{m}_k$ converges to $\Tilde{\mu} = (I,-\lambda D_pH(t,\cdot,D_x\Tilde{u},\Tilde{\mu})) \# \overline{m}$ uniformly in $[0,T]$. Since $\Tilde{u}_k \rightarrow \Tilde{u}$ in $C^{0,1}(Q)$, we get uniform bounds for $D_x\Tilde{u}_k$, and hence for $H(t,x,D_x\Tilde{u}_k,\Tilde{\mu}^k)$ and each of its first-order derivatives. Since
    $$
    \begin{aligned}
        &\frac{d}{dt}\int_\Omega |m - m_k|^2dx \\
        &\leq C_{\nu,\lambda}\left(\|D_pH\|_\infty^2\int_\Omega |m - m_k|^2 dx + \|m\|_{L^2(Q)}^2\|D_pH(t,x,D_x\Tilde{u},\Tilde{\mu}) - D_pH(t,x,D_x\Tilde{u}_k,\Tilde{\mu}^k)\|_\infty^2\right),
    \end{aligned}
    $$
    Gronwall's inequality gives
    $$
    \begin{aligned}
        \int_\Omega |m - m_k|^2dx &\leq C_{\nu,\lambda}T\|m\|_{L^2(Q)}^2e^{C_{\nu,\lambda}\|D_pH\|_\infty T}\|D_pH(t,x,D_x\Tilde{u},\Tilde{\mu}) - D_pH(t,x,D_x\Tilde{u}_k,\Tilde{\mu}^k)\|_\infty^2 \\
        &\leq C(\|D_x\Tilde{u} - D_x\Tilde{u}_k\|_\infty + \sup_{t \in [0,T]}d^*(\Tilde{\mu}(t),\Tilde{\mu}^k(t)))^2 \\
        &\rightarrow 0
    \end{aligned}
    $$
    and so $m_k \rightarrow m$ in $C^0(0,T;L^2(\Omega))$. Again, by Lemma \ref{Lem: Continuity of mu 1}, we get $\overline{\mu}^k \rightarrow \overline{\mu}$. Furthermore, by the results in Section \ref{Sec: Bootstrapping 1}, see in particular Lemma \ref{Lem: u C1+a}, we get uniform bounds for $\|u_k\|_{C^{1 +\alpha\beta^2/2,2+\alpha\beta^2}(Q)}$. By the Arzela-Ascoli theorem, there is a subsequence $(u_{k_j})_{j=1}^\infty$ converging to some $v$ in $C^{\alpha/2, 1 + \alpha}(Q)$. Since
    $$
    \begin{aligned}
        &2\nu\int_\Omega |D_xu - D_xu_k|^2dx - \frac{d}{dt}\int_\Omega |u - u_k|^2dx \\
        &\leq \int_\Omega |u - u_k|^2dx + \int_\Omega |H(t,x,D_x\Tilde{u},\mu) - H(t,x,D_x\Tilde{u}_k,\mu^k)|^2 dx,
    \end{aligned}
    $$
    it follows from Gronwall's inequality that
    $$
    \begin{aligned}
        \|u - u_k\|_{L^2(Q)}^2 \leq& Te^T(\|H(t,x,D_x\Tilde{u},\overline{\mu}) - H(t,x,D_x\Tilde{u}_k,\overline{\mu}^k)\|_{L^2(0,T;L^2(\Omega))}^2 \\
        &+ \|g(x,m(T)) - g(x,m_k(T))\|_{L^2(\Omega)}^2) \\
        \rightarrow& 0
    \end{aligned}
    $$
    and so $u_{k_j} \rightarrow u$ in $L^2(Q)$. By the uniqueness of the limit in $L^2(Q)$, this gives $v = u$ and so $(u_k)_{k=1}^\infty$ converges to $u$ in $C^{\alpha/2, 1 + \alpha}(Q)$.
    
    \textit{Compactness:} 
     Take $(\Tilde{u}_k,\Tilde{m}_k)_{k \in \N}$ bounded in $X$ and let $(u_k,m_k) = \Psi(\Tilde{u}_k,\Tilde{m}_k,\lambda)$ and $(u,m) = \Psi(\Tilde{u},\Tilde{m},\lambda)$. By similar arguments to those above, there is a subsequence $(u_{k_j})_{j=1}^\infty$ converging to some $u$ in $C^{0,1}(Q)$. Likewise, we get uniform bounds for $m_{k_j}$ in $C^{\alpha/2,\alpha}(Q)$, and so $m_{k_j}$ converges in $C^0(Q)$ (and hence $C^0(0,T;L^2(\Omega))$), passing to a subsequence if necessary.

    By the Leray-Schauder fixed point theorem, it follows that there exists some $(u,m) \in X$ with $\Psi(u,m,1) = (u,m)$. Letting
    $$\mu = (I,-D_pH(t,x,D_xu,\mu)) \# m,$$
    we get that $(u,m,\mu)$ is a classical solution to \eqref{Eq: MFGC}-\eqref{Eq: Dirchlet} (resp., \eqref{Eq: MFGC}-\eqref{Eq: Neumann}).
\end{proof}

\section{Dirichlet \& Neumann Boundary Conditions with Monotone Coupling} \label{sec: Problem 2}

For this problem, we will make most of our assumptions on the Lagrangian, and we will prove the existence of strong solutions to \eqref{Eq: MFGC} using an approach similar to that of \cite{kobeissi2022mean}, which leverages Lasry-Lions monotonicity to obtain a priori estimates. To this end, for $\theta \in (0,1]$, define $L^\theta(t,x,\alpha,\mu) = \theta L(t,x,\frac{\alpha}{\theta},\Theta(\mu))$ where $\Theta: \mathfrak{M}(\overline{\Omega} \times \R^n) \rightarrow \mathfrak{M}(\overline{\Omega} \times \R^n)$ is given by $\Theta(\mu) = (I \otimes \frac{1}{\theta}I) \# \mu$. Note that the associated Hamiltonian is given by $H^\theta(t,x,p,\mu) = \theta H(t,x,p,\Theta(\mu))$. Extending to $\theta = 0$ gives $H^0 = 0$ and
$$L^0(t,x,\alpha,\mu) = 
\begin{cases}
    0, &\alpha = 0 \\
    \infty, &\text{else}
\end{cases}$$
With the aim of applying the Leray-Schauder fixed point theorem, we will consider the following parametrized system:
\begin{equation}
    \begin{cases}
        -u_t - \nu\Delta u + H^\theta(t,x,D_xu,\mu) = \theta f(t,x,m), &(t,x) \in Q \\
        m_t - \nu\Delta m - \nabla \cdot (mD_pH^\theta(t,x,D_x u,\mu)) = 0, &(t,x) \in Q \\
        \mu = (I,-D_pH^\theta(t,\cdot,D_xu,\mu)) \# m, &t \in [0,T] \\
        u(T,x) = \theta g(x,m(T)), \hspace{1cm} m(0,x) = m_0(x), &x \in \overline{\Omega}
    \end{cases}
    \label{Eq: Parametrized MFGC 2}
\end{equation}
paired with either
\begin{equation}
    u = m = 0, \hspace{1cm} (t,x) \in \Sigma \tag{\ref{Eq: Parametrized MFGC 2}d}
    \label{Eq: Parametrized Dirchlet 2}
\end{equation}
or
\begin{equation}
    \frac{\partial u}{\partial\n} = \nu\frac{\partial m}{\partial\n} + mD_pH^\theta(t,x,D_xu,\mu) \cdot \n = 0, \hspace{1cm} (t,x) \in \Sigma \tag{\ref{Eq: Parametrized MFGC 2}n}
    \label{Eq: Parametrized Neumann 2}
\end{equation}

\subsection{Assumptions} \label{sec:assms2}

In the case of monotone couplings, we will still assume A\ref{A: Regularity of DpH}-\ref{A: Regularity of DxH}. However, we will replace A\ref{A: Hamiltonian}-\ref{A: DpH is Lipschitz in m} with the following assumptions:
\begin{assumption} \label{A: LL monotonicity}
    For all $R > 0$ and $\mu_1,\mu_2 \in \mathfrak{M}_{\infty,R}(\Omega \times \R^n)$, we have
    $$\int_{\Omega \times \R^n} (L(t,x,\alpha,\mu_1) - L(t,x,\alpha,\mu_2))d(\mu_1-\mu_2)(x,\alpha) \geq 0.$$
\end{assumption}
\begin{assumption} \label{A: Lagrangian}
    The Lagrangian $L: [0,T] \times \overline{\Omega} \times \R^n \times \mathfrak{M}(\overline{\Omega} \times \R^n) \rightarrow \R$ is differentiable with respect to $(x,\alpha)$, and $L$ and its derivatives are continuous on $[0,T] \times \overline{\Omega} \times \R^n \times \mathfrak{M}_{\infty,R}(\Omega \times \R^n)$ for any $R > 0$.
\end{assumption}
\begin{assumption} \label{A: Uniqueness of alpha}
    For every $t \in [0,T]$, $x \in \overline{\Omega}$, and $p \in \R^n$, the maximum in
    \begin{equation}
        \sup_{\alpha \in \R^n} \left[p \cdot \alpha - L(t,x,\alpha,\mu_\alpha)\right]
    \end{equation}
    is achieved at a unique $\alpha \in \R^n$.
\end{assumption}
\begin{assumption} \label{A: Lower bound for L}
    $L(t,x,\alpha,\mu) \geq \frac{1}{C_0}|\alpha|^{q'} - C_0\left(1 + \Lambda_{q'}(\mu)^{q'}\right)$ for some $q \in (1,\infty)$, $q' = \frac{q}{q-1}$, and $C_0 > 0$.
    In the Dirichlet case, we require $q \leq 2$.
\end{assumption}
\begin{assumption} \label{A: Bound for L and Lx}
    $|L(t,x,\alpha,\mu)| + |D_xL(t,x,\alpha,\mu)| \leq C_0\left(1 + |\alpha|^{q'} + \Lambda_{q'}(\mu)^{q'}\right)$.
\end{assumption}
\begin{assumption} \label{A: f, g, and m0}
    We take $m_0 \in C^{\beta}(\overline{\Omega})_+$ 
    with $\int_{\Omega} m_0 dx = 1$. For every $m \in \mathfrak{M}(\overline{\Omega})$, we have $g|_\Sigma = 0$ (resp., $\frac{\partial g}{\partial\n}|_\B = 0$) and $\|g(\cdot,m)\|_{C^{3+\beta}(\overline{\Omega})} + \|f(\cdot,\cdot,m)\|_{C^{\beta/2,1+\beta}(Q)} \leq C_0$.
    Furthermore, $m \mapsto g(\cdot,m)$ is continuous from $L^2(\Omega)$ into $L^2(\Omega)$, and $m \mapsto f(\cdot,\cdot,m)$ is continuous from $L^2(Q)$ into $L^2(Q)$.
\end{assumption}

As an example, one could consider a variation of the exhaustible resource model discussed in \cite{kobeissi2022mean} paired with Dirichlet boundary conditions. This would correspond to a situation in which players are forced to leave the game when their production capacities reach the boundary.

\begin{remark}
    Note that for $\theta \in (0,1]$, A\ref{A: Lower bound for L}-\ref{A: Bound for L and Lx} give
    $$L^\theta(t,x,\alpha,\mu) \geq \frac{\theta^{1-q'}}{C_0}|\alpha|^{q'} - C_0\theta - C_0\theta^{1-q'}\Lambda_{q'}(\mu)^{q'}$$
    $$|L^\theta(t,x,\alpha,\mu)| + |D_xL^\theta(t,x,\alpha,\mu)| \leq C_0\theta + C_0\theta^{1-q'}\left(1 + |\alpha|^{q'} + \Lambda_{q'}(\mu)^{q'}\right)$$
\end{remark}

\begin{remark}
    If we assume further regularity of $D_pH$ with respect to $\mu$ (e.g.~A\ref{A: DpH is Lipschitz in alpha}-\ref{A: DpH is Lipschitz in m}) such that we have H\"{o}lder estimates for $\mu$, we get classical solutions as in the previous problem. However, due to the nature of $\mu$ as a fixed-point, such estimates seem to require some kind of ``smallness condition". In this problem, we will instead show that we get strong (not classical) solutions under relatively modest assumptions.
\end{remark}

\subsection{Estimates on the Hamiltonian and Lagrangian}

We start by recalling results from \cite{kobeissi2022mean} for the Lagrangian and Hamiltonian, which rely heavily on properties of convex functions (see \cite{rockafellar1997convex}).

\begin{lemma} \label{Lem: Convexity of L}
    Fix $\theta \in (0,1]$. If $L^\theta$ is coercive and differentiable with respect to $\alpha$, then $L^\theta$ being strictly convex is equivalent to A\ref{A: Uniqueness of alpha}.
\end{lemma}

\begin{lemma} \label{Lem: Properties of H}
    Fix $\theta \in (0,1]$. Under assumptions A\ref{A: Lagrangian}-\ref{A: Bound for L and Lx}, the Hamiltonian $H^\theta(t,x,p,\mu)$ is differentiable with respect to $x$ and $p$, and $H$ and its derivatives are continuous on $[0,T] \times \overline{\Omega} \times \R^n \times \mathfrak{M}_{\infty,R}(\overline{\Omega} \times \R^n)$ for all $R > 0$. Furthermore, there exists $\Tilde{C}_0 > 0$ (depending only on $C_0$) so that
    \begin{equation}
        |D_pH^\theta(t,x,p,\mu)| \leq \Tilde{C}_0\theta(1 + |p|^{q-1} + \Lambda_{q'}(\mu))
        \label{Eq: Bound for DpH-theta}
    \end{equation}
    \begin{equation}
        |H^\theta(t,x,p,\mu)| \leq \Tilde{C}_0\theta(1 + |p|^{q}) + \Tilde{C}_0\theta^{1-q'}\Lambda_{q'}(\mu)^{q'}
        \label{Eq: Bound for H-theta}
    \end{equation}
    \begin{equation}
        p \cdot D_pH^\theta(t,x,p,\mu) - H^\theta(t,x,p,\mu) \geq \frac{\theta}{\Tilde{C}_0}|p|^q - \Tilde{C}_0\theta - \Tilde{C}_0\theta^{1-q'}\Lambda_{q'}(\mu)^{q'}
        \label{Eq: Convexity for H-theta}
    \end{equation}
    \begin{equation}
        |D_xH^\theta(t,x,p,\mu)| \leq \Tilde{C}_0\theta(1 + |p|^q) + \Tilde{C}_0\theta^{1-q'}\Lambda_{q'}(\mu)^{q'}
        \label{Eq: Bound for DxH-theta}
    \end{equation}
    for all $t \in [0,T]$, $x \in \overline{\Omega}$, $p \in \R^n$, and $\mu \in \mathfrak{M}(\overline{\Omega} \times \R^n)$. Without loss of generality, we will assume $\Tilde{C}_0 = C_0$.
\end{lemma}

\subsection{Fixed-Point Relation in $\mu$}
\label{Sec: Fixed-Point Relation 2}

As was done in \cite{kobeissi2022mean}, we will use our monotonicity assumption A\ref{A: LL monotonicity} to get a priori estimates for $\Lambda_\infty(\mu)$ and will prove the existence and uniqueness of the fixed-point $\mu$ using, respectively, the Leray-Schauder fixed point theorem and the monotonicity approach found in \cite[Lemma 5.2]{cardaliaguet2018mean}.

\begin{lemma} \label{Lem: Existence & Bound for mu 2}
    Assume A\ref{A: LL monotonicity}-\ref{A: Bound for L and Lx} hold. Given $t \in [0,T]$, $\theta \in (0,1]$, $p \in C^0(\Omega;\R^n)$, and $m \in \mathfrak{M}(\overline{\Omega})$, we have the following:
    \begin{enumerate}
        \item[1)] If $\mu$ satisfies
        \begin{equation}
            \mu = (I,-D_pH^\theta(t,\cdot,p(\cdot),\mu)) \# m,
            \label{Eq: mu 2}
        \end{equation}
        then we have
        \begin{equation}
            \Lambda_{q'}(\mu)^{q'} \leq 4C_0^2\theta^{q'} + \frac{(q')^{q-1}(2C_0)^q}{q}\theta^{q'}\|p\|_{L^q(m)}^q,
            \label{Eq: Bound for Lambda-q'}
        \end{equation}
        \begin{equation}
            \Lambda_\infty(\mu) \leq C_0\theta(1 + \|p\|_\infty + \Lambda_{q'}(\mu)).
            \label{Eq: Bound for Lambda-infinity}
        \end{equation}
        \item[2)] There exists a unique $\mu \in \mathfrak{M}_{\infty,R}(\Omega \times \R^n)$ satisfying \eqref{Eq: mu 2}, where $R \geq \|D_pH^\theta\|_\infty$.
    \end{enumerate}
\end{lemma}

\begin{proof}
    1) Applying A\ref{A: LL monotonicity} to $\mu$ and $m \otimes \delta_0$ gives
    $$\int_{\Omega \times \R^n} (L^\theta(t,x,\alpha,\mu) - L^\theta(t,x,\alpha,m \otimes \delta_0)) d\mu + \int_\Omega (L^\theta(t,x,0,m \otimes \delta_0) - L^\theta(t,x,0,\mu)) dm \geq 0.$$
    A\ref{A: Bound for L and Lx} gives
    $$\int_\Omega L^\theta(t,x,0,m \otimes \delta_0) dm \leq C_0\theta.$$
    Thus, by A\ref{A: Lower bound for L} and Lemma \ref{Lem: Convexity of L}, we get
    $$
    \begin{aligned}
        \frac{\theta^{1-q'}}{C_0}\int_{\Omega \times \R^n} |\alpha|^{q'}d\mu &\leq 2C_0\theta + \int_\Omega (L^\theta(t,x,\alpha^\mu,\mu) - L^\theta(t,x,0,\mu)) dm \\
        &\leq 2C_0\theta + \int_{\Omega \times \R^n} \alpha \cdot D_\alpha L^\theta(t,x,\alpha,\mu) d\mu \\
    \end{aligned}
    $$
    where $\mu = (I,\alpha^\mu) \# m$. Since $p(x) = -D_\alpha L^\theta(t,x,\alpha^\mu,\mu)$ and $yz \leq \frac{y^{q'}}{c^{q'}q'} + \frac{c^qz^q}{q}$ for $z,y \geq 0$ and $c > 0$, we get
    $$\frac{\theta^{1-q'}}{C_0}\int_{\Omega \times \R^n} |\alpha|^{q'}d\mu \leq 2C_0\theta + \frac{(2C_0q'\theta^{q'-1})^\frac{q}{q'}}{q}\int_\Omega |p(x)|^q dm + \frac{\theta^{1-q'}}{2C_0}\int_{\Omega \times \R^n} |\alpha|^{q'}d\mu.$$
    Since $\frac{q}{q'} + 1 = q$ and $(q' - 1)\frac{q}{q'} = 1$, this gives \eqref{Eq: Bound for Lambda-q'}. Combining this with \eqref{Eq: Bound for DpH-theta} gives \eqref{Eq: Bound for Lambda-infinity}.
    
    2) Take $\overline{p} \in C^0(\Omega;\R^n)$. For $\lambda \in [0,1]$, define
    $$L^{\overline{p},\lambda}(t,x,\alpha,\mu) = \lambda L^\theta(t,x,\alpha,\mu) + (1-\lambda)\left(\frac{|\alpha|^{q'}}{q'} - \alpha \cdot \overline{p}(x)\right)$$
    and denote by $H^{\overline{p},\lambda}$ the associated Legendre transform. As $|\alpha \cdot \overline{p}(x)| \leq \frac{|\alpha|^{q'}}{2q'} + \frac{2^{q-1}}{q}\|\overline{p}\|_\infty$, we can assume (up to a change of $C_0$) that $L^{\overline{p},\lambda}$ satisfies A\ref{A: LL monotonicity}-\ref{A: Bound for L and Lx}. Thus, $(\lambda,t,x,p,\mu) \mapsto D_pH^{\overline{p},\lambda}(t,x,p,\mu)$ is continuous on $[0,1] \times [0,T] \times \overline{\Omega} \times \R^n \times \mathfrak{M}_{\infty,R}(\Omega \times \R^n)$ for all $R > 0$.

    Define $\Phi: C^0(\overline{\Omega};\R^n) \times [0,1] \rightarrow C^0(\overline{\Omega};\R^n)$ by $\Phi(\alpha,\lambda) = -D_pH^{\overline{p},\lambda}(t,x,\overline{p},(I,\alpha) \# m)$. Note that $\Phi(\cdot,0)$ is constant, and that $\Phi$ is continuous by the continuity of the map $(\lambda,t,x,p,\mu) \mapsto D_pH^{\overline{p},\lambda}(t,x,p,\mu)$. For every $R > 0$, $A_R \coloneqq [0,1] \times [0,T] \times \overline{\Omega} \times B_R(0) \times \mathfrak{M}_{\infty,R}(\Omega \times \R^n)$ is compact, and hence $(\lambda,t,x,p,\mu) \mapsto D_pH^{\overline{p},\lambda}(t,x,p,\mu)$ is uniformly continuous on $A_R$ (by the Heine-Cantor theorem). Furthermore, $\overline{p}$ is uniformly continuous on $\overline{\Omega}$. Thus, the Arzela-Ascoli theorem gives that $\Phi$ is compact. Finally, if $\Phi(\alpha,\lambda) = \alpha$ for some $\lambda \in [0,1]$, then \eqref{Eq: Bound for Lambda-q'},\eqref{Eq: Bound for Lambda-infinity} give a uniform bound for $\|\alpha\|_\infty$. By the Leray-Schauder fixed point theorem, there exists some $\alpha_0$ such that $\Phi(\alpha_0,1) = \alpha_0$ and hence $\mu = (I,\alpha_0) \# m$ satisfies \eqref{Eq: mu 2}.

    To prove uniqueness, suppose $\mu_1,\mu_2$ both satisfy \eqref{Eq: mu 2}. Letting $\alpha_i = -D_pH^\theta(t,\cdot,p,\mu_i)$, A\ref{A: LL monotonicity} gives
    $$
    \begin{aligned}
        0 &\leq \int_{\Omega \times \R^n} (L^\theta(t,x,\alpha,\mu_1) - L^\theta(t,x,\alpha,\mu_2))d(\mu_1-\mu_2)(x,\alpha) \\
        &= \int_\Omega (L(t,x,\alpha_1,\mu_1) - L^\theta(t,x,\alpha_2,\mu_1) + L^\theta(t,x,\alpha_2,\mu_2) - L^\theta(t,x,\alpha_1,\mu_2)) dm \\
        &\leq \int_\Omega (D_\alpha L^\theta(t,x,\alpha_1,\mu_1) \cdot (\alpha_2 - \alpha_1) + D_\alpha L^\theta(t,x,\alpha_2,\mu_2) \cdot (\alpha_1 - \alpha_2)) dm.
    \end{aligned}
    $$
    As $D_\alpha L^\theta(t,x,\alpha_i,\mu_i) = p$, this shows that $\alpha_1 = \alpha_2$ $m$-a.e.
\end{proof}

We conclude this section by introducing two continuity results, the first will be useful in this section and the second will be useful in Section \ref{sec: Problem 3}. The proofs are omitted as they are nearly identical to that of Lemma \ref{Lem: Continuity of mu 1}.

\begin{lemma} \label{Lem: Continuity of mu 2}
    Assume A\ref{A: LL monotonicity}-\ref{A: Bound for L and Lx} hold. Fix $\theta \in (0,1]$ and let $(t_k,p_k,m_k)_{n \in \N}$ be a sequence in $[0,T] \times C^0(\overline{\Omega},\R^n) \times L^1(\Omega)$ such that
    \begin{itemize}
        \item $t_k \rightarrow t$ in $[0,T]$;
        \item $p_k \rightarrow p$ in $C^0(\overline{\Omega};\R^n)$;
        \item $\|m_k\|_{L^1(\Omega)} \leq 1$ and $m_k \rightarrow m$ in $L^1(\Omega)$.
    \end{itemize}
    Then $\mu^k \rightarrow \mu$ in $\mathfrak{M}_{\infty,R}(\Omega \times \R^n)$, where $\mu^k$ and $\mu$ are the fixed points associated to $(t_k,p_k,m_k)$ and $(t,p,m)$, respectively.
\end{lemma}

\begin{lemma} \label{Lem: New Continuity of mu 3}
    Assume A\ref{A: LL monotonicity}-\ref{A: Bound for L and Lx} hold. Fix $\theta \in (0,1]$ and let $(t_k,p_k,m_k)_{k \in \N}$ be a sequence in $[0,T] \times C^0(\Omega,\R^n) \times L^1(\Omega)$ such that
    \begin{itemize}
        \item $t_k \rightarrow t$ in $[0,T]$;
        \item $\abs{p_k} \leq C$ and $p_k \rightarrow p$ pointwise in $\Omega$;
        \item $\|m_k\|_{L^1(\Omega)} \leq C$ and $m_k \rightarrow m$ in $L^1(\Omega)$.
    \end{itemize}
    Then $\mu^k \rightarrow \mu$ in $\mathfrak{M}_{\infty,R}(\Omega \times \R^n)$, where $\mu^k$ and $\mu$ are the fixed points associated to $(t_k,p_k,m_k)$ and $(t,p,m)$, respectively.
\end{lemma}

\begin{remark} \label{rmk: uniform convergence of mu 2}
    This shows that for any $\theta \in (0,1]$, $u \in C^{0,1}([0,T] \times \Omega)$, and $m \in C^0(0,T;L^1(\Omega))$, the map $t \mapsto \mu(t) = (I,-D_pH^\theta(t,\cdot,D_xu(t,\cdot),\mu(t)) \# m(t)$ is continuous. Furthermore, if $(p_k,m_k)_{k \in \N}$ converges to $(D_xu,m)$ in $C^0(Q; \R^n) \times C^0(0,T;L^1(\Omega))$, then $\mu_k(t)$ defined by $$\mu_k(t) \coloneqq (I,-D_pH^\theta(t,\cdot,p_k(t,\cdot),\mu_k(t))) \# m_k(t)$$ converges to $\mu$ uniformly on $[0,T]$.
\end{remark}

\subsection{A Priori Estimates}
\label{Sec: A Priori Estimates}

In this section, we prove a priori bounds on $|u|$, $|D_xu|$, and $\Lambda_\infty(\mu)$ using an approach similar to that found in \cite{kobeissi2022mean}, which we adapt to account for the boundary conditions. This approach leverages our Lasry-Lions monotonicity assumption, allowing us to eliminate some of the ``smallness conditions" that we required in the non-monotone case. As before, we start with the Neumann case.

\begin{theorem} \label{Thm: A Priori Estimates (N)}
    Assume A\ref{A: LL monotonicity}-\ref{A: f, g, and m0} hold and suppose $(u,m,\mu)$ is a strong solution to \eqref{Eq: Parametrized MFGC 2}-\eqref{Eq: Parametrized Neumann 2}. Then there exists some $C > 0$ (depending only on the constants in the assumptions) such that $\|u\|_\infty \leq C\theta$. Furthermore, if $u \in C^{1,3}([0,T] \times \Omega)$, then up to a change of constants, $\|D_xu\|_\infty \leq C\theta^\frac{1}{2}$ and $\underset{t \in [0,T]}{\sup} \Lambda_\infty(\mu(t)) \leq C\theta$.
\end{theorem}

\begin{proof}
    \textit{Estimating $\int_0^T\Lambda_{q'}(\mu)^{q'}dt$: }
    
    Define $(X,\alpha)$ by
    $$\begin{cases}
        \alpha_t = -D_pH^\theta(t,X_t,D_xu(t,X_t),\mu(t)) \\
        dX_t = \alpha_tdt + \sqrt{2\nu}dB_t - K_t \\
        X_0 = \xi \sim m_0
    \end{cases}$$
    where $\{B_t\}_{t \in [0,T]}$ is a Brownian motion independent of $\xi$ and $K_t$ enforces the reflection boundary conditions. Note that for $t \in [0,T]$ and $s \in [t,T]$,
    $$\alpha_s = \underset{\alpha'}{\text{argmin}} E\left[\int_t^T (L^\theta(s,X_s^{\alpha'},\alpha'_s,\mu(s)) + \theta f(x,X_s^{\alpha'},m(s)))ds + \theta g(X_T^{\alpha'},m(T))\right]$$
    where
    $$\begin{cases}
        dX_t^{\alpha'} = \alpha'_tdt + \sqrt{2\nu}dB'_t - K_t^{\alpha'} \\
        X_0^{\alpha'} = \xi' \sim m_0
    \end{cases}$$
    $\{B'_t\}_{t \in [0,T]}$ is a Brownian motion independent of $\xi'$, and $K_t^{\alpha'}$ enforces the reflection boundary conditions. Define $\Tilde{X}$ by
    $$\begin{cases}
        d\Tilde{X}_t = \sqrt{2\nu}dB_t - \Tilde{K}_t \\
        \Tilde{X}_0 = \xi \sim m_0
    \end{cases}
    $$
    and define the measures $\Tilde{m}(t) = \mathcal{L}(\Tilde{X}_t)$ and $ \Tilde{\mu}(t) = \mathcal{L}(\Tilde{X}_t) \otimes \delta_0$, where $\mathcal{L}(\Tilde{X}_t)$ denotes the law of $\Tilde{X}_t$.

    Taking $\alpha' = 0$ gives
    $$
    \begin{aligned}
        &\int_0^T\int_{\Omega \times \R^n} L^\theta(t,x,\alpha,\mu(t)) d\mu(t,x,\alpha)dt + \int_0^T\int_\Omega \theta f(t,x,m(t)) dm(t,x)dt + \int_\Omega \theta g(x,m(T)) dm(T,x) \\
        &\leq \int_0^T\int_{\Omega \times \R^n} L^\theta(t,x,\alpha,\mu(t)) d\Tilde{\mu}(t,x,\alpha)dt + \int_0^T\int_\Omega \theta f(t,x,m(t)) d\Tilde{m}(t,x)dt + \int_\Omega \theta g(x,m(T)) d\Tilde{m}(T,x).
    \end{aligned}
    $$
    By A\ref{A: f, g, and m0},
    $$\int_0^T\int_{\Omega \times \R^n} L^\theta(t,x,\alpha,\mu(t)) d\mu(t,x,\alpha)dt \leq \int_0^T\int_{\Omega \times \R^n} L^\theta(t,x,\alpha,\mu(t)) d\Tilde{\mu}(t,x,\alpha)dt + 2C_0\theta(1 + T)$$
    Furthermore, by A\ref{A: LL monotonicity} and A\ref{A: Bound for L and Lx}, we get
    $$
    \begin{aligned}
        &\int_{\Omega \times \R^n} L^\theta(t,x,\alpha,\mu(t)) d\Tilde{\mu}(t,x,\alpha) + \int_{\Omega \times \R^n} L^\theta(t,x,\alpha,\Tilde{\mu}(t)) d\mu(t,x,\alpha) \\
        &\leq \int_{\Omega \times \R^n} L^\theta(t,x,\alpha,\mu(t)) d\mu(t,x,\alpha) + \int_{\Omega \times \R^n} L^\theta(t,x,\alpha,\Tilde{\mu}(t)) d\Tilde{\mu}(t,x,\alpha)
    \end{aligned}
    $$
    and
    $$\int_{\Omega \times \R^n} L^\theta(t,x,\alpha,\Tilde{\mu}(t)) d\Tilde{\mu}(t,x,\alpha) = \theta\int_\Omega L(t,x,0,\Tilde{\mu}(t)) \Tilde{m}(t,dx) \leq C_0\theta.$$
    By A\ref{A: Lower bound for L}, these estimates give
    $$
    \begin{aligned}
        \int_0^T \Lambda_{q'}(\mu)^{q'} dt &= \int_0^T\int_{\Omega \times \R^n} |\alpha|^{q'} d\mu(t,x,\alpha)dt \\
        &\leq C_0^2\theta^{q'}T + C_0\theta^{q'-1}\int_0^T\int_{\Omega \times \R^n} L^\theta(t,x,\alpha,\Tilde{\mu}(t)) d\mu(t,x,\alpha)dt \\
        &\leq 2C_0^2\theta^{q'}(1 + 2T).
    \end{aligned}
    $$
    
    \textit{Estimating $\|u\|_\infty$: }
    For $w(t,x) = u(t,x) + 2C_0\theta(t-T) - C_0\theta^{1-q'}\int_t^T \Lambda_{q'}(\mu) dt$, we have
    $$
    \begin{aligned}
        \frac{d}{dt}\int_\Omega (w-k)_+^2 dx &\geq 2\int_\Omega \left(\nu|D_xw|^2 + (w-k)\int_0^1 D_pH^\theta(t,x,sD_xu,\mu(t)) ds \cdot D_xw\right)\chi_{w \geq k}dx \\
        &\geq -C_\nu\|D_pH^\theta\|_\infty^2\int_\Omega (w-k)_+^2 dx.        
    \end{aligned}
    $$
    By Gronwall's inequality,
    $$\int_\Omega (w(t,x)-k)_+^2 dx \leq e^{C_\nu\|D_pH^\theta\|_\infty^2(T-t)}\int_\Omega (w(T,x)-k)_+^2 dx.$$
    Choosing $k = \underset{x}{\max} \, w(T,x)_+ = \underset{x}{\max} \, u(T,x)_+$ gives $u(t,x) \leq \underset{x}{\max} \, u(T,x)_+ + 2C_0\theta(T - t) + C_0\theta^{1-q'}\int_t^T \Lambda_{q'}(\mu) dt$. Similarly, $u(t,x) \geq -\underset{x}{\max} \, u(T,x)_- + 2C_0\theta(t-T) - C_0\theta^{1-q'}\int_t^T \Lambda_{q'}(\mu) dt$. In particular,
    $$\|u\|_\infty \leq C_0\theta(1 + 2T) + C_0\theta^{1-q'}\int_0^T \Lambda_{q'}(\mu) dt.$$

    \textit{Estimating $\|D_xu\|_\infty$: }
    Note that if $u \in C^{1,3}([0,T] \times \Omega)$, then
    $$
    \begin{aligned}
        &-\frac{1}{2}\partial_t|D_xu|^2 - \nu\Delta D_xu \cdot D_xu + D_{xx}^2u D_pH^\theta(t,x,D_xu,\mu) \cdot D_xu + D_xH^{\theta}(x,D_xu,\mu) \cdot D_xu = D_xf \cdot D_xu.
    \end{aligned}
    $$
    If $\|u\|_\infty = 0$, then $u \equiv 0$. So suppose $\|u\|_\infty > 0$. Define the following functions:
    $$\varphi(v) = \exp\left(\exp\left(-1-R^{-1} v\right)\right), \text{ for } |v| \leq R \coloneqq \|u\|_\infty + (1 + C_0)\theta^{1-q'}\int_0^T\Lambda_{q'}(\mu)^{q'}dt,$$
    $$w(t,x) = \varphi\left(u(T-t,x) + (1 + C_0)\theta^{1-q'}\int_0^t\Lambda_{q'}(\mu(s))^{q'}ds\right)|D_xu(T-t,x)|^2.$$
    
    Recall that
    $$1 \leq \varphi(v) \leq e, \hspace{.7cm} bR^{-1}e^{-2} \leq \frac{|\varphi'(v)|}{\varphi(v)} \leq bR^{-1}, \hspace{.7cm} \frac{|\varphi'|}{\varphi^{1+\frac{q}{2}}} \geq bR^{-1}e^{-2}e^{-\frac{q}{2}}.$$
    Furthermore, by similar arguments to those in the proof of Theorem \ref{Thm: Gradient Estimate (N)},
    $$
    \begin{aligned}
        &\partial_tw - \nu\Delta w + D_xw \cdot D_pH^\theta(t,x,D_xu,\mu) - 2\nu\frac{\varphi'}{\varphi}D_xw \cdot D_xu + 2\nu\varphi|D_{xx}^2 u|^2 \\
        &= \frac{\varphi'}{\varphi}w\left[\theta f + D_xu \cdot D_pH^\theta(t,x,D_xu,\mu) - H^\theta(t,x,D_xu,\mu) + (1 + C_0)\theta^{1-q'}\Lambda_{q'}(\mu)^{q'}\right] \\
        &\hspace{0.5cm}- \nu\frac{\varphi''\varphi-2(\varphi')^2}{\varphi^3}w^2 - 2\varphi D_xu \cdot D_xH^{\theta} + 2\theta\varphi D_xf \cdot D_xu.
    \end{aligned}
    $$
    We now bound the right-hand side from above. To begin, notice that $\varphi''\varphi - 2(\varphi')^2 \geq 0$ for $a \geq b$. By \eqref{Eq: Bound for DxH-theta},
    $$-2\varphi D_xu \cdot D_xH^{\theta} \leq 2C_0\varphi|D_xu|\left(\theta + \theta|D_xu|^q + \theta^{1-q'}\Lambda_{q'}(\mu)^{q'}\right).$$
    Since $\varphi' < 0$, \eqref{Eq: Convexity for H-theta} gives
    $$
    \begin{aligned}
        &\frac{\varphi'}{\varphi}w\left[\theta f + D_xu \cdot D_pH^\theta(t,x,D_xu,\mu) - H^\theta(t,x,D_xu,\mu) + (1 + C_0)\theta^{1-q'}\Lambda_{q'}(\mu)^{q'}\right] \\
        &\leq -C_0^{-1}\theta\frac{\varphi'}{\varphi^{1+\frac{q}{2}}}w^{1+\frac{q}{2}} + 2C_0\theta\frac{|\varphi'|}{\varphi}w - \theta^{1-q'}\frac{|\varphi'|}{\varphi}w\Lambda_{q'}(\mu)^{q'}.
    \end{aligned}
    $$
    Combining these inequalities, we get
    \begin{equation}
        \begin{aligned}
            &\partial_tw - \nu\Delta w + D_xw \cdot D_pH(t,x,D_xu,\mu) - 2\nu\frac{\varphi'}{\varphi}D_xw \cdot D_xu + 2\nu\varphi|D_{xx}^2 u|^2 \\
            &\leq \theta C_0R^{-1}\bigg[-\frac{e^{-2}e^{-\frac{q}{2}}}{C_0^2}w^{1+\frac{q}{2}} + 2w - \frac{1}{C_0\theta^{q'}}\Lambda_{q'}(\mu)^{q'}w \\
            &\qquad+ 2Rw^\frac{1}{2}\left(2e^{\frac{1}{2}} + w^\frac{q}{2} + e^{\frac{1}{2}}\theta^{-q'}\Lambda_{q'}(\mu)^{q'}\right)\bigg] \\
            &\leq \theta C_0R^{-1}\bigg[\left(\varepsilon - \frac{e^{-2}e^{-\frac{q}{2}}}{C_0^2}\right)w^{1+\frac{q}{2}} + C_{\varepsilon,q}\left(1 + 2R\right)  \\
            &\qquad+ \left(\varepsilon - \frac{1}{C_0\theta^{q'}}\right)\Lambda_{q'}(\mu)^{q'}w + \left(C_{\varepsilon} + 2R\right)\theta^{-q'}\Lambda_{q'}(\mu)^{q'}\bigg] \\
            &\leq \theta \Tilde{C}_{\varepsilon,q}\left(\theta^{-q'}\Lambda_{q'}(\mu)^{q'} + 1\right)
        \end{aligned}
        \label{Eq: w inequality 2}
    \end{equation}
    for all $\varepsilon > 0$. Note that $v = w - \Tilde{C}_{\varepsilon,q}\left(\theta^{1-q'}\int_0^t \Lambda_{q'}(\mu(s))^{q'}ds + \theta t\right)$ satisfies
    $$\partial_tv - \nu \Delta v + D_xv \cdot \left(D_pH^\theta - 2\nu\frac{\varphi'}{\varphi}D_xu\right) \leq 0.$$
    
    Recall that for all $z \in C^2(\overline{\Omega})$ with $\frac{\partial z}{\partial\n} = 0$ on $\B$, we have $\frac{\partial}{\partial\n}|D_xz|^2 \leq 0$ on $\B$. Hence, $\frac{\partial}{\partial\n}w \leq 0$. Since
    $$\frac{d}{dt}\int_\Omega (v-k)_+^2 dx + 2\nu\int_\Omega |D_xv|^2\chi_{v \geq k} dx \leq 2\int_\Omega (v-k)_+D_xv \cdot \left(D_pH^\theta - 2\nu\frac{\varphi'}{\varphi}D_xu\right) dx,$$
    we have
    $$\frac{d}{dt}\int_\Omega (v-k)_+^2 dx \leq C_\nu\left(\|D_pH^\theta\|_\infty^2 + 4\nu^2\left\|\frac{\varphi'}{\varphi}\right\|_\infty^2\|D_xu\|_\infty^2\right)\int_\Omega (v-k)_+^2 dx.$$
    By Gronwall's inequality, it follows that $v \leq \|v(0,\cdot)\|_\infty$ and so
    $$w \leq \|w(0,\cdot)\|_\infty + \Tilde{C}_{\varepsilon,q}\left(\theta^{1-q'}\int_0^T \Lambda_{q'}(\mu)^{q'}dt + \theta T\right) \leq \theta C.$$
    Since $\varphi \geq 1$, this shows $\|D_xu\|_\infty \leq \theta^\frac{1}{2}C$.

    \textit{Estimating $\Lambda_\infty(\mu)$: }
    Combining our gradient estimate with \eqref{Eq: Bound for Lambda-q'},\eqref{Eq: Bound for Lambda-infinity} gives $\Lambda_{q'}(\mu) + \Lambda_\infty(\mu) \leq C\theta$ for some constant $C$ depending only on the constants in the assumptions.
\end{proof}

\begin{theorem} \label{Thm: A Priori Estimates (D)}
    Assume A\ref{A: LL monotonicity}-\ref{A: f, g, and m0} hold and suppose $(u,m,\mu)$ is a strong solution to \eqref{Eq: Parametrized MFGC 2}-\eqref{Eq: Parametrized Dirchlet 2}. Then there exists some $C > 0$ (depending only on the constants in the assumptions) such that $\|u\|_\infty \leq C\theta$. Furthermore, if $u \in C^{1,3}([0,T] \times \Omega)$, then $\|D_xu\|_\infty \leq C$ and $\underset{t \in [0,T]}{\sup} \Lambda_\infty(\mu(t)) \leq C\theta$.
\end{theorem}

\begin{proof}
    The proof is nearly identical to the proof of Theorem \ref{Thm: A Priori Estimates (N)}, with an estimate for the gradient on the boundary given by an argument similar to the proof of Lemma \ref{Lem: Bounds for Du on Boundary (D)}.
\end{proof}

\subsection{Bootstrapping}
\label{Sec: Bootstrapping 2}

Let $X \coloneqq C^{0,1}(Q) \times C^0(0,T;L^2(\Omega))$. With the aim of applying the Leray-Schauder fixed point theorem, we will define the map $\Gamma: X \times [0,1] \rightarrow X$ so that $(u,m) = \Gamma(\Tilde{u},\Tilde{m},\theta)$ is the unique strong solution to
    \begin{equation}
        \begin{cases}
            -u_t - \nu\Delta u + H^\theta(t,x,D_x\Tilde{u},\mu) = \theta f(t,x,m), &(t,x) \in Q \\
            m_t - \nu\Delta m - \nabla \cdot (mD_pH^\theta(t,x,D_x\Tilde{u},\mu)) = 0, &(t,x) \in Q \\
            \mu = (I,-D_pH^\theta(t,\cdot,D_x\Tilde{u},\mu)) \# \overline{m}, &t \in [0,T] \\
            u(T,x) = \theta g(x,m(T)), \hspace{1cm} m(0,x) = m_0(x), &x \in \overline{\Omega}
        \end{cases}
        \label{Eq: Definition of Gamma}
    \end{equation}
    paired with either
\begin{equation}
    u = m = 0, \hspace{1cm} (t,x) \in \Sigma \tag{\ref{Eq: Definition of Gamma}d}
    \label{Eq: Gamma Dirchlet}
\end{equation}
or
\begin{equation}
    \frac{\partial u}{\partial\n} = \nu\frac{\partial m}{\partial\n} + mD_pH^\theta(t,x,D_x\Tilde{u},\mu) \cdot \n = 0, \hspace{1cm} (t,x) \in \Sigma \tag{\ref{Eq: Definition of Gamma}n}
    \label{Eq: Gamma Neumann}
\end{equation}
(where $\overline{m}$ is defined as in section \ref{Sec: Bootstrapping 1}).

\begin{remark} \label{rmk: Holder Estimates for u}
    As in section \ref{Sec: Bootstrapping 1}, $m$ and $\mu$ are well-defined. Moreover, \cite[Theorem IV.9.1]{ladyzhenskaia1968linear} (resp., the discussion at the end of \cite[Section IV.9]{ladyzhenskaia1968linear}) gives that there is a unique strong solution $u$ in $W^{1,2}_{\Tilde{q}}(Q)$ for arbitrarily large $\Tilde{q}$. Hence, $u \in C^{\frac{\beta}{2},1+\beta}(Q)$ with bounds depending only on $\|D_x\Tilde{u}\|_\infty$. Thus, $\Gamma$ is well-defined. 
\end{remark}

\begin{lemma}
    Assume A\ref{A: LL monotonicity}-\ref{A: f, g, and m0} hold. Then there exist $\alpha \in (0,\beta)$ and $C > 0$ so that $m \in C^{\alpha/2,\alpha}(Q)$ with
    $$\|m\|_{C^{\alpha/2,\alpha}(Q)} \leq C,$$
    where $C$ depends only on $\|D_x\Tilde{u}\|_\infty$ and the constants in the assumptions.
\end{lemma}

\begin{theorem} \label{Lem: Schauder Estimates for Strong Solutions}
    Assume A\ref{A: Regularity of DpH}-\ref{A: f, g, and m0} hold. If $(u,m,\mu)$ is a strong solution to \eqref{Eq: Parametrized MFGC 2}-\eqref{Eq: Parametrized Dirchlet 2} (resp., \eqref{Eq: Parametrized MFGC 2}-\eqref{Eq: Parametrized Neumann 2}), then we have
    $$\|u\|_{C^{\beta/2,1+\beta}(Q)} \leq C,$$
    where $C$ depends only on the constants in the assumptions.
\end{theorem}

\begin{proof}
    Let $\alpha_k$ be a sequence in $C^{\beta_k/2,\beta_k}(Q;\R^n)$ converging to $-D_pH^\theta(t,x,D_xu,\mu)$ with $\|\alpha_k\|_\infty \leq C$. Letting $\theta \in [0,1]$ and $\mu^k = (I,\alpha_k) \# m$, arguments similar to those in section \ref{sec: Problem 1} give $C^{1,3}([0,T] \times \Omega)$ classical solutions to the HJ equation
    \begin{equation}
        \begin{cases}
            -\partial_tu_k - \nu\Delta u_k + H^\theta(t,x,D_xu_k,\mu^k) = \theta f(t,x,m), &(t,x) \in Q \\
            u_k(T,x) = \theta g(x,m(T)), &x \in \overline{\Omega}
        \end{cases}
        \label{Eq: Sequence of HJs}
    \end{equation}
    paired with $u_k\vert_\Sigma = 0$ (resp., $\frac{\partial u_k}{\partial\n}\vert_\Sigma = 0$). By nearly identical arguments to those used to prove our a priori estimates, this gives uniform bounds for $\|D_xu_k\|_\infty$ and hence $\|u_k\|_{C^{\beta/2,1+\beta}(Q)}$. By the Arzela-Ascoli theorem, we get that there is a subsequence $u_{k_j}$ converging to some $v \in C^{0,1}(Q)$.

    Since
    $$
    \begin{aligned}
        &2\nu\int_\Omega |D_xu_{k_j} - D_xu|^2 dx - \frac{d}{dt}\int_\Omega |u - u_{k_j}|^2 dx \\
        &= 2\int_\Omega (u - u_{k_j})(H^\theta(t,x,D_xu_{k_j},\mu^{k_j}) - H^\theta(t,x,D_xu,\mu)) dx \\
        &\leq \nu\int_\Omega |D_xu_{k_j} - D_xu|^2 dx + C_\nu(1+\|D_pH^\theta\|_\infty^2)\int_\Omega |u - u_{k_j}|^2 dx + Cd^*(\mu(t),\mu^{k_j}(t))^{2\beta}
    \end{aligned}
    $$
    for a.e. $0 \leq t \leq T$, Gronwall's inequality gives that
    $$
    \begin{aligned}
        \int_\Omega |u - u_{k_j}|^2 dx &\leq C\int_t^Td^*(\mu(\tau),\mu^{k_j}(\tau))^{2\beta} d\tau \\
        &\leq C\int_0^T\left(\int_\Omega |D_pH^\theta(t,x,D_xu,\mu) + \alpha_{k_j}|dm\right)^{2\beta}d\tau  \\
        &\rightarrow 0
    \end{aligned}
    $$
    for a.e. $0 \leq t \leq T$. By uniqueness, it follows that $u = v$ and so $\|u\|_{C^{\beta/2,1+\beta}(Q)} \leq C$ (by Remark \ref{rmk: Holder Estimates for u}).
\end{proof}

\subsection{Existence}

\begin{theorem}
    Under assumptions A\ref{A: Regularity of DpH}-\ref{A: f, g, and m0}, there is a strong solution to \eqref{Eq: MFGC}-\eqref{Eq: Dirchlet} (resp., \eqref{Eq: MFGC}-\eqref{Eq: Neumann}).
\end{theorem}

\begin{proof}
    Define $X$ and $\Gamma: X \times [0,1] \rightarrow X$ as in section \ref{Sec: Bootstrapping 2}.

    \textit{$\Gamma(\cdot,0)$ is constant:} 
    First, note that for every $(\Tilde{u},\Tilde{m}) \in X$, $\Gamma(\Tilde{u},\Tilde{m},0) = (0,\rho)$, where $\rho$ is the weak solution to the heat equation $\rho_t = \nu\Delta\rho$ with $\rho(0,x) = m_0$ and $\rho|_\Sigma = 0$ (resp., $\frac{\partial \rho}{\partial\n}|_\Sigma = 0$).
    
    \textit{Bound for fixed-points:} 
    By the results from the previous sections, we get
    $$\|u\|_{C^{\beta/2,1+\beta}(Q)} + \|m\|_{C^0(0,T;L^2(\Omega))} \leq C$$
    for all $(u,m,\theta) \in X \times [0,1]$ with $\Gamma(u,m,\theta) = (u,m)$.

    \textit{Continuity:} Fix $\theta \in [0,1]$ and take $(\Tilde{u}_k,\Tilde{m}_k) \rightarrow (\Tilde{u},\Tilde{m})$ in $X$. Define $(u_k,m_k) = \Gamma(\Tilde{u}_k,\Tilde{m}_k,\theta)$ and $(u,m) = \Gamma(\Tilde{u},\Tilde{m},\theta)$. By Remark \ref{rmk: uniform convergence of mu 2}, it follows that $\mu^k = (I,-D_pH^\theta(t,\cdot,D_x\Tilde{u}_k,\mu^k)) \# \overline{m}_k$ converges to $\mu = (I,-D_pH^\theta(t,\cdot,D_x\Tilde{u},\mu)) \# \overline{m}$ uniformly in $[0,T]$. Since $D_x\Tilde{u}_k$ is uniformly bounded, so are $H^\theta(t,x,D_x\Tilde{u}_k,\mu^k)$ and each of its first-order derivatives.
    
    By similar arguments to those in the proof of Theorem \ref{Thm: Existence 1}, we get $m_k \rightarrow m$ in $C^0(0,T;L^2(\Omega))$. Furthermore, by the results in Section \ref{Sec: Bootstrapping 2}, see in particular Lemma \ref{Lem: Schauder Estimates for Strong Solutions}, we get uniform bounds for $\|u_k\|_{C^{\beta/2, 1 + \beta}(Q)}$. By the Arzela-Ascoli theorem, there is a subsequence $(u_{k_j})_{j=1}^\infty$ converging to some $v$ in $C^{0,1}(Q)$.  By similar arguments to those in the proof of Theorem \ref{Thm: Existence 1}, $v = u$ and so $(u_k)_{k=1}^\infty$ converges to $u$ in $C^{0,1}(Q)$.
    
    \textit{Compactness:} 
     Take $(\Tilde{u}_k,\Tilde{m}_k)_{k \in \N}$ bounded in $X$ and let $(u_k,m_k) = \Gamma(\Tilde{u}_k,\Tilde{m}_k,\theta)$ and $(u,m) = \Gamma(\Tilde{u},\Tilde{m},\theta)$. By similar arguments to those above, there is a subsequence $(u_{k_j})_{j=1}^\infty$ converging to some $u$ in $C^{0,1}(Q)$. Likewise, we get uniform bounds for $m_{k_j}$ in $C^{\alpha/2,\alpha}(Q)$, and so $m_{k_j}$ converges in $C^0(Q)$ (and hence $C^0(0,T;L^2(\Omega))$), passing to a subsequence if necessary.

    By the Leray-Schauder fixed point theorem, it follows that there exists some $(u,m) \in X$ with $\Gamma(u,m,1) = (u,m)$. Letting
    $$\mu = (I,-D_pH(t,x,D_xu,\mu)) \# m,$$
    we get that $(u,m,\mu)$ is a strong solution to \eqref{Eq: MFGC}-\eqref{Eq: Dirchlet} (resp., \eqref{Eq: MFGC}-\eqref{Eq: Neumann}).
\end{proof}

\section{Uniqueness}
\label{Sec: Uniqueness}

In this section, we prove uniqueness of solutions to our boundary-value problems under two different types of assumptions. For the first case, we use the approach found in \cite{kobeissi2022mean}, which relies on the monotonicity assumptions A\ref{A: LL monotonicity} and
\begin{assumption} \label{A: Monotonicity for g,f}
    For all $t \in [0,T]$ and $m_1,m_2 \in \mathfrak{M}(\overline{\Omega})$, we have
    $$\int_\Omega (g(x,m_1) - g(x,m_2))(m_1-m_2) dx \geq 0$$
    and
    $$\int_\Omega (f(t,x,m_1) - f(t,x,m_2))d(m_1-m_2)(x) \geq 0.$$
\end{assumption}

\begin{theorem}
    Assume A\ref{A: LL monotonicity},A\ref{A: Monotonicity for g,f} and either A\ref{A: Hamiltonian} or A\ref{A: Lagrangian}-\ref{A: Uniqueness of alpha} hold. Then there is at most one strong solution $(u,m,\mu)$ to \eqref{Eq: MFGC}-\eqref{Eq: Dirchlet} (resp., \eqref{Eq: MFGC}-\eqref{Eq: Neumann}).
\end{theorem}

\begin{proof}
    Let $(u_1,m_1,\mu_1)$ and $(u_2,m_2,\mu_2)$ be strong solutions. By A\ref{A: Monotonicity for g,f}, we get
    $$
    \begin{aligned}
        0 \leq& \int_\Omega (g(x,m_1(T)) - g(x,m_2(T)))(m_1(T,x) - m_2(T,x)) dx \\
        &+ \int_0^T\int_\Omega (f(t,x,m_1) - f(t,x,m_2))(m_1-m_2) dxdt \\
        =& \int_0^T\int_\Omega m_1(H(t,x,D_xu_1,\mu_1) - H(t,x,D_xu_2,\mu_2) + D_x(u_2 - u_1) \cdot D_pH(t,x,D_xu_1,\mu_1)) dxdt \\
        &+ \int_0^T\int_\Omega m_2(H(t,x,D_xu_2,\mu_2) - H(t,x,D_xu_1,\mu_1) + D_x(u_1 - u_2) \cdot D_pH(t,x,D_xu_2,\mu_2)) dxdt.
    \end{aligned}
    $$
    Note that for $i = 1,2$,
    $$L(t,x,\alpha^{\mu_i},\mu_i) = D_xu_i \cdot D_pH(t,x,D_xu_i,\mu_i) - H(t,x,D_xu_i,\mu_i)$$
    and
    $$D_xu_i = -D_\alpha L(t,x,\alpha^{\mu_i},\mu_i).$$
    Thus,
    $$
    \begin{aligned}
        0 \leq& \int_0^T\int_\Omega m_1(L(t,x,\alpha^{\mu_2},\mu_2) - L(t,x,\alpha^{\mu_1},\mu_1) + D_\alpha L(t,x,\alpha^{\mu_2},\mu_2) \cdot (\alpha^{\mu_1}-\alpha^{\mu_2})) dxdt \\
        &+ \int_0^T\int_\Omega m_2(L(t,x,\alpha^{\mu_1},\mu_1) - L(t,x,\alpha^{\mu_2},\mu_2) + D_\alpha L(t,x,\alpha^{\mu_1},\mu_1) \cdot (\alpha^{\mu_2}-\alpha^{\mu_1})) dxdt.
    \end{aligned}
    $$
    Since $L$ is strictly convex,
    \begin{equation}
        L(t,x,\alpha_1,\mu) - L(t,x,\alpha_2,\mu) + D_\alpha L(t,x,\alpha_1,\mu) \cdot (\alpha_2 - \alpha_1) \leq 0
        \label{Eq: Convexity of L}
    \end{equation}
    with equality holding if and only if $\alpha_1 = \alpha_2$. Hence,
    $$
    \begin{aligned}
        0 \leq& \int_0^T\int_\Omega \bigg(m_1(L(t,x,\alpha^{\mu_1},\mu_2) - L(t,x,\alpha^{\mu_1},\mu_1)) + m_2(L(t,x,\alpha^{\mu_2},\mu_1) - L(t,x,\alpha^{\mu_2},\mu_2))\bigg)dxdt \\
        =& -\int_0^T\int_{\Omega \times \R^n} (L(t,x,\alpha,\mu_1) - L(t,x,\alpha,\mu_2)) d(\mu_1 - \mu_2)(x,\alpha)dt
    \end{aligned}
    $$
    By A\ref{A: LL monotonicity}, this gives
    $$\int_0^T\int_\Omega \bigg(m_1(L(t,x,\alpha^{\mu_2},\mu_1) - L(t,x,\alpha^{\mu_1},\mu_1)) + m_2(L(t,x,\alpha^{\mu_1},\mu_2) - L(t,x,\alpha^{\mu_2},\mu_2))\bigg)dxdt = 0.$$
    By the condition for equality for \eqref{Eq: Convexity of L}, we get $|\{(t,x) \in Q : \alpha^{\mu_1} \neq \alpha^{\mu_2}, m_i \neq 0\}| = 0$ for $i = 1,2$. Therefore, $\alpha^{\mu_1} = \alpha^{\mu_2}$. By the uniqueness of solutions to the Fokker-Planck equation, $m_1 = m_2$. Therefore, $\mu_1 = (I,\alpha^{\mu_i}) \# m_i = \mu_2$. By the uniqueness of solutions to the Hamilton-Jacobi equation, $u_1 = u_2$.
\end{proof}

For the second case, we refrain from assuming Lasry-Lions monotonicity and instead use an approach adapted from \cite{kobeissi2022classical} to prove the uniqueness of classical solutions. This requires a short time horizon but may be more realistic for some models (e.g.~models of crowd dynamics).

\begin{theorem}
    Suppose that $D_pH$ is Lipschitz in $p$ and that we have uniform bounds on $\|u\|_{C^{1,2}(Q)}$ and $\|m\|_{C^{\alpha,\alpha/2}(Q)}$ for solutions $(u,m,\mu)$ to \eqref{Eq: MFGC}-\eqref{Eq: Dirchlet} (resp., \eqref{Eq: MFGC}-\eqref{Eq: Neumann}) (as is the case if A\ref{A: Hamiltonian}-\ref{A: Regularity of DxH} hold).
    Then there exists $T_0 > 0$ such that if $T \leq T_0$ and
    \begin{enumerate}
        \item For all $m_1,m_2 \in \mathfrak{M}(\overline{\Omega})$, we have
        $$\|g(\cdot,m_1) - g(\cdot,m_2)\|_{C^{1+\beta}(\Omega)} \leq C_0d^*(m_1,m_2);$$
        \item For every $p \in \R^n$, $\alpha_1,\alpha_2 \in C^0(\overline{\Omega};\R^n)$, and $m_1,m_2 \in \mathfrak{M}(\overline{\Omega})$, we have
        $$|H(t,x,p,(I,\alpha_1) \# m_1) - H(t,x,p,(I,\alpha_2) \# m_2)| \leq C_0(\|\alpha_1 - \alpha_2\|_\infty + d^*(m_1,m_2)),$$
        $$|D_pH(t,x,p,(I,\alpha_1) \# m_1) - D_pH(t,x,p,(I,\alpha_2) \# m_1)| \leq L_1\|\alpha_1 - \alpha_2\|_\infty,$$
        $$|D_pH(t,x,p,(I,\alpha_1) \# m_1) - D_pH(t,x,p,(I,\alpha_1) \# m_2)| \leq C_0d^*(m_1,m_2)$$
        for some $L_1 \in (0,1)$;
    \end{enumerate}
     then there is at most one classical solution  to \eqref{Eq: MFGC}-\eqref{Eq: Dirchlet} (resp., \eqref{Eq: MFGC}-\eqref{Eq: Neumann}).
\end{theorem}

\begin{proof}
    Let $(u_1,m_1,\mu_1)$ and $(u_2,m_2,\mu_2)$ be solutions. Then there exist $\alpha \in (0,1)$ and $C > 0$ (depending only on the constants in the assumptions) so that
    $$\|u_1\|_{C^{1,2}(Q)} + \|u_2\|_{C^{1,2}(Q)} + \|m_1\|_{C^{\alpha/2,\alpha}(Q)} + \|m_2\|_{C^{\alpha/2,\alpha}(Q)} \leq C.$$
    Define $u \coloneqq u_1 - u_2$, $m \coloneqq m_1 - m_2$, $\alpha_i(t,x) \coloneqq -D_pH(t,x,D_xu_i(t,x),\mu_i)$, and $\alpha \coloneqq \alpha_1 - \alpha_2$. Then
    $$
    \begin{aligned}
        |\alpha| \leq& |D_pH(t,x,D_xu_1,\mu_1) - D_pH(t,x,D_xu_1,(I,\alpha_1) \# m_2)| \\
        &+ |D_pH(t,x,D_xu_1,(I,\alpha_1) \# m_2) - D_pH(t,x,D_xu_1,\mu_2)| \\
        &+ |D_pH(t,x,D_xu_1,\mu_2) - D_pH(t,x,D_xu_2,\mu_2)| \\
        \leq& C_0d^*(m_1,m_2) + L_1|\alpha| + C|D_xu|
    \end{aligned}
    $$
    and
    $$
    \begin{aligned}
        d^*(m_1(t),m_2(t)) &\leq \max\{1,diam(\Omega)\}\int_\Omega |m_1(t) - m_2(t)|dx \\
        &\leq \max\{1,diam(\Omega)\}|\Omega|^\frac{1}{2}\sup_{0 \leq s \leq T}\|m_1(s) - m_2(s)\|_{L^2(\Omega)} \\
        &\leq C_\nu T\|\alpha\|_\infty
    \end{aligned}
    $$
    and so $\|\alpha\|_\infty \leq C\|D_xu\|_\infty$ provided $T$ is sufficiently small. Now note that $u$ solves
    $$\begin{cases}
        -u_t - \nu\Delta u + H(t,x,D_xu_1,\mu_1) - H(t,x,D_xu_2,\mu_2) = 0, &(t,x) \in Q \\
        u(T,\cdot) = g(\cdot,m_1(T)) - g(\cdot,m_2(T)), &x \in \overline{\Omega}
    \end{cases}$$
    (with either Dirichlet or Neumann boundary conditions). Note that
    $$|H(t,x,D_xu_1,\mu_1) - H(t,x,D_xu_2,\mu_2)| \leq C\sup_{0 \leq t \leq T} \|D_xu(t)\|_\infty.$$
    Furthermore, by Theorem 6.48 (resp., Theorem 6.49) in \cite{lieberman1996second}, since $u(T,\cdot) \in C^{1 + \beta}(\overline{\Omega})$, we get that $u \in C^{\frac{1+\beta}{2},1+\beta}(Q)$ with
    $$\|D_xu\|_\infty \leq \|D_xu(T,\cdot)\|_\infty + CT^\frac{\beta}{2}(\|D_xu\|_\infty + \|u(T,\cdot)\|_{C^{1+\beta}(\Omega)}).$$
    Combining this with our assumptions gives
    $$
    \begin{aligned}
        \|D_xu\|_\infty &\leq C_0d^*(m_1,m_2) + CT^{\beta/2}(\|D_xu\|_\infty + C_0d^*(m_1,m_2)) \\
        &\leq C(T + T^{\beta/2} + T^{1 + \beta/2})\|D_xu\|_\infty
    \end{aligned}
    $$
    Thus, choosing $T$ sufficiently small gives $D_xu = 0$ and hence $m = 0$. This implies that $m_1 = m_2$, $\mu_1 = \mu_2$, and finally $u_1 = u_2$.
\end{proof}

\section{Invariance Constraints} \label{sec: Problem 3}

In this final section, we shift our attention to MFGCs under invariance constraints imposed on the state space.

\subsection{Assumptions} \label{sec: assms 3}

In the case of invariance constraints, we will still assume A\ref{A: Regularity of DpH}-\ref{A: Bound for L and Lx}, and A\ref{A: Monotonicity for g,f} hold. However, we will replace the rest of our assumptions with the following:
\begin{assumption} \label{A: a}
    There exists a sequence of non-negative numbers $\lambda_k$ such that
    \begin{enumerate}
        \item $a$ is uniformly bounded in $C^{1+\beta}(\overline{\Omega})^{n \times n}$;
        \item $\sum_{i,j} (a(x))_{ij} \xi_i\xi_j \geq \lambda_k|\xi|^2$ for every $\xi \in \R^n$ and $x \in \Omega_{1/k}$.
    \end{enumerate}
\end{assumption}
\begin{assumption} \label{A: sigma}
    There exists a matrix $\sigma \in W^{1,\infty}(\Omega)$ such that $a(x) = \sigma(x)\sigma^*(x)$.
\end{assumption}
\begin{assumption} \label{A: m0 and f}
    $m_0 \in C^{\beta/2,\beta}(\overline{\Omega})$ with $m_0 \geq 0$ and $\int_{\Omega} m_0 dx = 1$. For every $m \in C^0(0,T;\mathfrak{M}(\Omega))$, we have $\|f(\cdot,\cdot,m)\|_{C^{\beta/2,1+\beta}(Q)} \leq C_0$. Furthermore, the map $m \mapsto f(\cdot,\cdot,m)$ is continuous from $L^2(Q)$ into $L^2(Q)$.
\end{assumption}
\begin{assumption} \label{A: g}
     The map $g: L^1(\Omega) \rightarrow C^0(\Omega)$ given by $m \mapsto g(\cdot,m)$ is such that $\frac{\partial g}{\partial\n}|_\B = 0$ and the map $m \mapsto g(\cdot,m)$ is continuous from $L^2(\Omega)$ into $L^2(\Omega)$. Furthermore, there exist a sequence $\{\Omega^\varepsilon\}$ of convex, $C^{2+\alpha}$ subdomains such that $\Omega_{2\varepsilon} \subseteq \Omega^\varepsilon \subseteq \Omega_\varepsilon$ and a sequence of functions $g_\varepsilon(x,m)$ satisfying
    \begin{enumerate}
        \item For every $\varepsilon > 0$ and $m \in \mathfrak{M}(\Omega^\varepsilon)$, we have $\frac{\partial g_\varepsilon}{\partial\n_\varepsilon}|_{\partial\Omega^\varepsilon} = 0$, $\|g_\varepsilon(\cdot,m)\|_{C^{3+\beta}(\overline{\Omega^\varepsilon})} \leq C(\varepsilon)$, and $\|g_\varepsilon(\cdot,m)\|_{C^{1+\beta}(\overline{\Omega^\varepsilon})} \leq C_0$;
        \item For every $\varepsilon > 0$, $m \mapsto g_\varepsilon(\cdot,m)$ is continuous from $L^2(\Omega^\varepsilon)$ into $L^2(\Omega^\varepsilon)$;
        \item If $m_\varepsilon \rightarrow m$ in $L^1(\Omega)$, then (up to a subsequence) $g_\varepsilon(\cdot,m_\varepsilon) \rightarrow g(\cdot,m)$ a.e.;
    \end{enumerate}
\end{assumption}
For an example satisfying Assumption \ref{A: g}, let $\B$ be $C^{4+\beta}$-smooth, $m \mapsto g(\cdot,m)$ be continuous from $L^1(\Omega)$ into $L^2(\Omega)$, and $\|g(\cdot,m)\|_{C^{3+\beta}(\overline{\Omega})} \leq C$ for all $m$:
    Then there exists a collar region $\{x: \abs{d(x)} < r\}$ on which $d(x)$ is $C^{4+\beta}$.
    Let $\varepsilon < r/2$.
    Pick a $C^\infty$ smooth non-increasing function $\eta_\varepsilon:\intco{0,\infty} \to \intco{0,\infty}$ such that $\eta_\varepsilon(s) = \varepsilon$ for $s \in [0,\varepsilon]$ and $\eta_\varepsilon$ has compact support in $\intco{0,2\varepsilon}$.
    Define $\xi = \xi(t,x)$ to be the flow map given by solving the ODE
    \begin{equation} \label{Eq: Definition of xi}
        \od{\xi(t)}{t} = \eta_\varepsilon\del{d\del{\xi(t)}}Dd\del{\xi(t)}, \quad \xi(0) = x.
    \end{equation}
    As the vector field $V(\xi) = \eta_\varepsilon\del{d(\xi)}Dd(\xi)$ is $C^{3+\beta}$ smooth, the map $x \mapsto \xi(t,x)$ is a $C^{3+\beta}$ diffeomorphism for every $t \geq 0$.
    Since $\abs{Dd} = 1$, we have
    \begin{equation}
        \od{}{t}d\del{\xi(t)} = Dd\del{\xi(t)}\cdot \od{\xi(t)}{t} = \eta_\varepsilon\del{d\del{\xi(t)}}.
    \end{equation}
    From the structure of $\eta_\varepsilon$ we can deduce that $d\del{\xi(t,x)} \geq \varepsilon$ for all $t \geq 1$ and all $x \in \overline{\Omega}$, and on the other hand $\xi(t,x) = x$ whenever $d(x) \geq 2\varepsilon$.
    We now set $\psi_\varepsilon(x) = \xi(1,x)$.
    Then $\psi_\varepsilon$ is a $C^{3+\beta}$ diffeomorphism.
    Setting $\Omega^\varepsilon \coloneqq \psi_\varepsilon(\Omega)$,
    we have $\Omega_{2\varepsilon} \subseteq \Omega^\varepsilon \subseteq \Omega_\varepsilon$.
    If we set $g_\varepsilon(\cdot,m) \coloneqq g(\cdot,m) \circ \psi_\varepsilon^{-1}$, we have the desired properties.
\begin{assumption} \label{A: Invariance Condition}
    There exist some $\delta > 0$ and $C > 0$ so that $H$ satisfies \eqref{Eq: Invariance Condition} for all $p \in \R^n$ and all $(t,x) \in [0,T] \times \Gamma_\delta$. Without loss of generality, we may assume that $\delta < \varepsilon_0$, where $\varepsilon_0$ is the constant given in Definition \ref{Def: Subdomains}.
\end{assumption}

\begin{remark} \label{Rmk: Neumann Problem}
    In the previous sections, we proved the well-posedness of MFGCs with Neumann boundary conditions in the case where $a = \nu I$. However, the arguments can be generalized easily enough to the case of variable coefficients $a_{ij} \in C^{1+\beta}(\overline{\Omega})$, provided they satisfy the uniform ellipticity constraint $a_{ij}(x)\xi_i \xi_j \geq \nu|\xi|^2$ for some constant $\nu > 0$. On the other hand, in Assumption A\ref{A: a} we take the ellipticity constraint to degenerate as $x$ approaches the boundary of $\Omega$. The reason for this is as follows. Consider the approximating system without variable coefficients in the diffusion:
    \begin{equation*}
        \begin{cases}
            -\partial_tu_\varepsilon - \nu\Delta u_\varepsilon + H(t,x,D_xu_\varepsilon,\mu_\varepsilon) = f(t,x,m_\varepsilon), &(t,x) \in [0,T] \times \Omega^\varepsilon \\
            \partial_tm_\varepsilon - \nu \Delta m_\varepsilon - \nabla \cdot (m_\varepsilon D_pH(t,x,D_x u_\varepsilon,\mu_\varepsilon)) = 0, &(t,x) \in [0,T] \times \Omega^\varepsilon \\
            \mu_\varepsilon(t) = (I,-D_pH(t,\cdot,D_xu_\varepsilon,\mu_\varepsilon)) \# m(t,\cdot), &t \in [0,T] \\
            u_\varepsilon(T,x) = g_\varepsilon(x,m_\varepsilon(T)), \hspace{1cm} m_\varepsilon(0,x) = m_0(x), &x \in \Omega^\varepsilon \\
            D_xu_\varepsilon \cdot \n_\varepsilon = 0, \hspace{1cm} [\nu D_xm_\varepsilon + m_\varepsilon D_pH(t,x,D_xu_\varepsilon,\mu_\varepsilon)] \cdot \n_\varepsilon = 0, &(t,x) \in [0,T] \times \partial \Omega^\varepsilon.
        \end{cases}
    \end{equation*}
    By careful review of the proofs in the previous sections, we get a solution $(u_\varepsilon,m_\varepsilon,\mu_\varepsilon)$ with uniform bounds on $D_xu_\varepsilon$ and hence on $D_pH(t,x,D_xu_\varepsilon,\mu_\varepsilon)$.
    As $\varepsilon \to 0$ this is incompatible with the invariance constraint \eqref{Eq: Invariance Condition} (see \cite[Remark 4.4]{porretta2020mean}).
\end{remark}

\subsection{Example (Cournot)}
\label{Sec: Example}

To motivate the application of invariance constraints to MFGCs, we will consider a modified version of the Cournot mean field game system (see \cite{camilli2025learning}) without a discount term, which we adjust to allow for degeneracy of the diffusion coefficient. Let $\Omega = (0,L)$ for some $L > 0$. Now let $a(x) = \sigma(x)^2$ for some $\sigma \in C^2([0,L])$ such that
\begin{enumerate}
    \item $\sigma(0) = \sigma(L) = 0$ and $\sigma(x) \geq \nu(\varepsilon) > 0$ in $(\varepsilon,L-\varepsilon)$ for every $\varepsilon > 0$.
    \item There exist $C_\sigma,\delta > 0$ so that $\sigma(x) \leq C_\sigma x$ in $(0,\delta)$ and $\sigma(x) \leq C_\sigma(L-x)$ in $(L-\delta,L)$.
\end{enumerate}
Let $c \in C^2([0,L])$ be a function such that $c(x) > 0$ for $x > 0$ and $c(x) \leq x$ in $[0,\varepsilon]$ for some $\varepsilon > 0$.

Consider the system
\begin{equation} \label{Eq: Example}
    \begin{cases}
        -\partial_tu - a\partial_{xx}^2u + \underset{v \geq 0}{\sup} \{v\partial_xu + v\pi(t) - \gamma v - \frac{v^2}{c(x)}\} = 0, &\text{in } [0,T] \times (0,L) \\
        \partial_tm - \partial_{xx}^2(am)  - \partial_x(mq^*) = 0, &\text{in } [0,T] \times (0,L) \\
        \pi(t) = P\left(t,\int_0^Lq^*(t,x)m(t,x) dx\right), &\text{in } [0,T] \\
        u(T,x) = g(x), \qquad m(0,x) = m_0(x), &\text{in } (0,L) \\
        q^*(t,x) = \underset{v \geq 0}{\operatorname{argmax}}\{v\partial_xu + v\pi(t) - \gamma v - \frac{v^2}{c(x)}\}
    \end{cases}
\end{equation}
where $P(t,\cdot)$ is a decreasing function for all $t \in [0,T]$.

Note that
$$L(t,x,\widetilde{\alpha},(I,\alpha) \# \rho) =
\begin{cases}
    -\widetilde{\alpha} P\left(t,\int_0^L \alpha\rho dx\right) + \gamma\widetilde{\alpha} + \frac{\widetilde{\alpha}^2}{c(x)}, &\widetilde{\alpha} \geq 0 \\
    \infty, &\widetilde{\alpha} < 0
\end{cases}$$
satisfies A\ref{A: LL monotonicity} for $\mu_1,\mu_2 \in \mathfrak{M}([0,L] \times [0,\infty))$, as
$$
\begin{aligned}
    &\int_{[0,L] \times [0,\infty]}(L(t,x,\alpha,\mu_1) - L(t,x,\alpha,\mu_2))d(\mu_1-\mu_2)(x,\alpha) \\
    &= \left(P\left(t,\int_0^L \alpha_2\rho_2 dx\right) - P\left(t,\int_0^L \alpha_1\rho_1 dx\right)\right)\int_0^L (\alpha_1\rho_1 - \alpha_2\rho_2) dx \\
    &\geq 0
\end{aligned}
$$
where $\mu_i = (I,\alpha_i) \# \rho_i$.
Furthermore, the Hamiltonian
$$H(t,x,p,(I,\alpha) \# \rho) = \underset{v \geq 0}{\sup} \left\{vp + vP\left(t,\int_0^L \alpha\rho dx\right) - \gamma v - \frac{v^2}{c(x)}\right\}$$
satisfies
$$D_pH(t,x,p,\mu) = q^*(t,x) = \left(\frac{c(x)(p+P+\gamma)}{2}\right)_+$$
and
$$D_xH(t,x,p,\mu) = c'(x)\left[\frac{q^*(t,x)}{c(x)}\right]^2 = c'(x)\left(\frac{p+P+\gamma}{2}\right)_+^2$$
a.e.
In Lemma \ref{Lem: Gradient Bound} we will establish bounds on $D_x u$ and $\mu$ that are independent of the diffusion coefficients.
Therefore, to establish the invariance condition \eqref{Eq: Invariance Condition}, it is enough to show that there is some constant $C > 0$ such that
$\sigma(x)^2 \leq -xD_pH(t,x,p,\mu) + Cx^2$ in $(0,\delta)$ and $\sigma(x)^2 \leq (L-x)D_pH(t,x,p,\mu) + C(L-x)^2$ in $(L-\delta,L)$ for all $p,\mu$ satisfying $|p|, |P| \leq M$, where $M$ is defined a priori and $P$ is defined in terms of $\mu$ as above.
Indeed, on $(L-\delta,L)$ we just need $C \geq C_\sigma^2$, while in $(0,\delta)$ we observe that $-xD_pH(t,x,p,\mu) = -xc(x)\left(\frac{(p+P+\gamma)}{2}\right)_+ \geq -x^2(M + \gamma/2)$, so it suffices to take $C = C_\sigma^2 + M + \gamma/2$.    

We note that the assumptions A\ref{A: Regularity of DpH}-\ref{A: Lagrangian} and A\ref{A: Invariance Condition} are not technically satisfied (e.g.~we need to restrict ourselves to non-negative controls and $p,\mu$ satisfying given a priori estimates). However, the well-posedness of this system follows by nearly identical arguments to those found in the following sections.

\subsection{Existence of Solutions}
\label{Sec: Existence 3}

In this section, we will prove the following existence result:

\begin{theorem} \label{Thm: Existence}
    Assume A\ref{A: Regularity of DpH}-\ref{A: Bound for L and Lx} and A\ref{A: Monotonicity for g,f}-\ref{A: Invariance Condition} hold. Then there is a weak solution to \eqref{Eq: MFGC 2}.
\end{theorem}

Throughout this section, we will assume A\ref{A: Regularity of DpH}-\ref{A: Bound for L and Lx} and A\ref{A: Monotonicity for g,f}-\ref{A: Invariance Condition} hold. For $\varepsilon > 0$, define $(u_\varepsilon, m_\varepsilon, \mu_\varepsilon)$ to be the unique strong solution to
\begin{equation}
    \begin{cases}
        -\partial_tu_\varepsilon - \sum_{i,j} a_{ij}\partial_{ij}^2 u_\varepsilon + H(t,x,D_xu_\varepsilon,\mu_\varepsilon) = f(t,x,m_\varepsilon), &(t,x) \in [0,T] \times \Omega^\varepsilon \\
        \partial_tm_\varepsilon - \sum_{i,j} \partial_{ij}^2 (a_{ij}m_\varepsilon) - \nabla \cdot (m_\varepsilon D_pH(t,x,D_x u_\varepsilon,\mu_\varepsilon)) = 0, &(t,x) \in [0,T] \times \Omega^\varepsilon \\
        \mu_\varepsilon(t) = (I,-D_pH(t,\cdot,D_xu_\varepsilon,\mu_\varepsilon)) \# m_\varepsilon(t,\cdot), &t \in [0,T] \\
        u_\varepsilon(T,x) = g_\varepsilon(x,m_\varepsilon(T)), \hspace{1cm} m_\varepsilon(0,x) = m_0(x), &x \in \Omega^\varepsilon \\
        aD_xu_\varepsilon \cdot \n_\varepsilon = 0, \hspace{1cm} [a^*D_xm_\varepsilon + m_\varepsilon (\Tilde{b}(x) + D_pH(t,x,D_xu_\varepsilon,\mu_\varepsilon))] \cdot \n_\varepsilon = 0, &(t,x) \in [0,T] \times \partial \Omega^\varepsilon
    \end{cases}
    \label{Eq: Approximating MFGC}
\end{equation}
We extend to $Q$ by setting $u_\varepsilon = m_\varepsilon = 0$ on $Q \setminus [0,T] \times \Omega^\varepsilon$. Note that $m_\varepsilon \geq 0$ with $\int_{\Omega^\varepsilon} m_\varepsilon dx = \int_{\Omega^\varepsilon} m_0 dx \leq 1$ for all $t \in [0,T]$. By A\ref{A: m0 and f} and A\ref{A: g}, $f(t,x,m_\varepsilon)$ and $g_\varepsilon(x,m_\varepsilon)$ are uniformly bounded.

By arguments nearly identical to those used to prove Theorem \ref{Thm: A Priori Estimates (N)}, we get
\begin{equation} \label{Eq: A Priori Estimates}
    \|u_\varepsilon\|_\infty + \int_0^T \Lambda_{q'}(\mu_\varepsilon)^{q'} dt \leq C
\end{equation}
and
\begin{equation} \label{Eq: Gradient Bound Depending on Epsilon}
    \|D_xu_\varepsilon\|_\infty + \Lambda_{q'}(\mu_\varepsilon)^{q'} + \Lambda_\infty(\mu_\varepsilon) \leq C_\varepsilon.
\end{equation}
Furthermore, adapting the argument in \cite{porretta2020mean} for Lipschitz regularity, we can actually bound $C_\varepsilon$ uniformly in $\varepsilon$.

\begin{lemma} \label{Lem: Gradient Bound}
    Assume A\ref{A: Regularity of DpH}-\ref{A: Bound for L and Lx} and A\ref{A: Monotonicity for g,f}-\ref{A: Invariance Condition} hold. Then there exists some constant $C > 0$ independent of $\varepsilon$ such that
    $$\|D_xu_\varepsilon\|_\infty \leq C.$$
\end{lemma}

\begin{proof}
    We will assume that $u_\varepsilon$ is smooth (otherwise, we can approximate $u_\varepsilon$ as in Lemma \ref{Lem: Schauder Estimates for Strong Solutions}). As was done in \cite{porretta2020mean}, we will define $w_\varepsilon = |D_xu_\varepsilon|^2e^{d^\gamma}$ for some $\gamma \in (0,1)$. Then we have
    $$
    \begin{aligned}
        -\partial_tw_\varepsilon - tr(aD_{xx}^2w_\varepsilon) =& 2e^{d^\gamma}D_xu_\varepsilon \cdot (D_xf - D_xH - D_pHD_{xx}^2u_\varepsilon) - 2\gamma d^{\gamma - 1}(aD_xd) \cdot D_xw_\varepsilon \\
        &+ \left(\gamma^2d^{2\gamma-2} - \gamma(\gamma - 1)d^{\gamma - 2}\right)w_\varepsilon(aD_xd) \cdot D_xd - \gamma d^{\gamma-1}tr(aD_{xx}^2d)w_\varepsilon \\
        &+ 2e^{d^\gamma}\sum_{i,j,k} (\partial_ka_{ij}\partial_ku_\varepsilon\partial_{ij}^2u_\varepsilon - a_{ij}\partial_{jk}^2u_\varepsilon\partial_{ik}^2u_\varepsilon).
    \end{aligned}
    $$
    Assumption A\ref{A: m0 and f} and \eqref{Eq: Gradient Bound Depending on Epsilon} give
    $$D_xu_\varepsilon \cdot (D_xf - D_xH) \leq C_\varepsilon$$
    and assumption A\ref{A: sigma} gives
    $$
    \begin{aligned}
        \sum_{i,j,k} (\partial_ka_{ij}\partial_ku_\varepsilon\partial_{ij}^2u_\varepsilon - a_{ij}\partial_{jk}^2u_\varepsilon\partial_{ik}^2u_\varepsilon) &= 2\sum_{i,j,k,l}\sigma_{il}\partial_k\sigma_{il}\partial_ku_\varepsilon\partial_{ij}^2u_\varepsilon - \sum_k |\sigma^*D_x\partial_ku_\varepsilon|^2 \\
        &\leq C|D_xu_\varepsilon|^2 \\
        &\leq C_\varepsilon.
    \end{aligned}
    $$
    Thus,
    $$-\partial_tw_\varepsilon - tr(aD_{xx}^2w_\varepsilon) + (D_pH + 2\gamma d^{\gamma-1}(aD_xd)) \cdot D_xw_\varepsilon \leq c_\varepsilon w_\varepsilon + C_\varepsilon,$$
    where
    $$c_\varepsilon \coloneqq \gamma d^{\gamma-1}(D_pH \cdot D_xd - tr(aD_{xx}^2d)) + (\gamma^2d^{2\gamma-2} - \gamma(\gamma - 1) d^{\gamma-2})(aD_xd) \cdot D_xd.$$
    By A\ref{A: Invariance Condition}, we get
    $$
    \begin{aligned}
        c_\varepsilon &\leq - \gamma d^{\gamma-2}(aD_xd) \cdot D_xd + (\gamma^2d^{2\gamma-2} - \gamma(\gamma - 1) d^{\gamma-2})(aD_xd) \cdot D_xd + C_0d^\gamma + C \\
        &= \gamma^2(d-1)(aD_xd) \cdot D_xd + C_0d^\gamma + C \\
        &\leq C
    \end{aligned}
    $$
    in $\overline{\Omega} \setminus \Omega_\delta$. Combining this with \eqref{Eq: Gradient Bound Depending on Epsilon}, we get that
    \begin{equation} \label{Eq: Inequality for w}
        -\partial_tw_\varepsilon - tr(a_\varepsilon D_{xx}^2w_\varepsilon) + (D_pH + 2\gamma d^{\gamma-1}(a_\varepsilon D_xd)) \cdot D_xw_\varepsilon \leq C_\varepsilon(w_\varepsilon + 1).
    \end{equation}
    By \cite[Lemma 4]{leonori2011gradient}, 
    the maximum of $w_\varepsilon$ must not be attained on $[0,T] \times \B_\varepsilon$. Furthermore, applying the maximum principle to \eqref{Eq: Inequality for w}, we get that the maximum must be attained on $\{T\} \times \Omega^\varepsilon$. Thus, A\ref{A: g} gives the desired estimate.
\end{proof}

With this, we are ready to prove our existence result.

\begin{proof}[Proof of Theorem \ref{Thm: Existence}]
    By Lemma \ref{Lem: Gradient Bound}, we get uniform bounds for $H$ and each of its first-order derivatives, which in turn gives uniform bounds for $u_\varepsilon$ in $C^{\beta/2,1+\beta}([0,T] \times K)$ for each $K \subset\subset \Omega$. By the Arzela-Ascoli theorem and a diagonal argument, there is some $u \in C^{0,1}([0,T] \times \Omega)$ such that $u_\varepsilon \rightarrow u$ in $C^{0,1}([0,T] \times K)$ for every $K \subset\subset \Omega$.
    
    Similarly, for every $\varepsilon_0 > 0$ and every $\varepsilon \in (0,\varepsilon_0)$, $m_\varepsilon$ is uniformly bounded in $C^{\alpha/2,1+\alpha}([0,T] \times \overline{\Omega}_{\varepsilon_0})$, and hence we have a subsequence converging to some $m$ uniformly on $[0,T] \times \overline{\Omega}_{\varepsilon_0}$. By a diagonal argument, there is a subsequence (which we will still denote by $m_\varepsilon$) converging to some $m$ in $C^{0,1}([0,T] \times K)$ for $K \subset\subset \Omega$.
    Note that $m \geq 0$ and $\int_\Omega m(t,\cdot) dx \leq 1$ for all $t \in [0,T]$.
    
    By nearly identical arguments to those used to prove \cite[Proposition 4.3]{porretta2020mean}, we get that for each $t \in [0,T]$,
    $$\int_\Omega m_\varepsilon(t,x) dx \rightarrow \int_\Omega m(t,x) dx$$
    and hence $m_\varepsilon(t,\cdot) \rightarrow m(t,\cdot)$ in $L^1(\Omega)$ (by Scheff\'{e}'s lemma).
    Thus, applying Lemma \ref{Lem: New Continuity of mu 3} gives that $\mu_\varepsilon(t)$ converges to
    $$\mu(t) = (I,-D_pH(t,\cdot,D_xu(t,\cdot),\mu(t))) \# m(t,\cdot)$$
    for all $t \in [0,T]$. By Lemma \ref{Lem: Properties of H}, this gives that $D_pH(t,x,D_xu_\varepsilon,\mu_\varepsilon) \rightarrow D_pH(t,x,D_xu,\mu)$ pointwise. By \cite[Proposition 4.3]{porretta2020mean}, $m_\varepsilon \rightarrow m$ in $C^0(0,T;L^1(\Omega))$ and $m$ is a weak solution to the FP equation. By A\ref{A: m0 and f} and A\ref{A: g}, $f(t,x,m_\varepsilon) \rightarrow f(t,x,m)$ a.e. and $g(x,m_\varepsilon(T)) \rightarrow g(x,m(T))$ in $L^p(\Omega)$ for $p < \infty$. By \cite[Proposition 3.4]{porretta2020mean}, $u$ is a weak solution to the HJ equation, thus completing the proof.
\end{proof}

Notice that as an immediate consequence of Lemma \ref{Lem: Gradient Bound}, we have a gradient estimate for the solutions constructed in the proof of Theorem \ref{Thm: Existence}. Moreover, by the uniqueness of solutions (see Section \ref{Sec: Uniqueness 3}), this gives an a priori bound for the gradient of solutions.

\subsection{Uniqueness of Solutions}
\label{Sec: Uniqueness 3}

Now that we have proven the existence of solutions, we shift our focus to uniqueness. The first half of our uniqueness proof will follow the argument used in \cite{porretta2020mean}. However, for the second half, we apply the strategy used in \cite{kobeissi2022mean} to adapt the argument to the MFGC setting.

\begin{theorem} \label{Thm: Uniqueness}
    Assume A\ref{A: Regularity of DpH}-\ref{A: Bound for L and Lx} and A\ref{A: Monotonicity for g,f}-\ref{A: Invariance Condition} hold. Then there is at most one solution to \eqref{Eq: MFGC 2}.
\end{theorem}

\begin{proof}
    Let $(u,m,\mu)$ and $(v,\rho,\nu)$ be solutions to \eqref{Eq: MFGC 2}. Defining $\Omega^\varepsilon$ as before, let $(u_\varepsilon,m_\varepsilon,\mu_\varepsilon)$ and $(v_\varepsilon,\rho_\varepsilon,\nu_\varepsilon)$ be solutions to the approximating systems
    \begin{equation}
        \begin{cases}
            -\partial_tu_\varepsilon - \sum_{i,j} a_{ij}\partial_{ij}^2 u_\varepsilon + H(t,x,D_xu_\varepsilon,\mu_\varepsilon) = f(t,x,m), &(t,x) \in [0,T] \times \Omega^\varepsilon \\
            \partial_tm_\varepsilon - \sum_{i,j} \partial_{ij}^2 (a_{ij}m_\varepsilon) - \nabla \cdot (m_\varepsilon D_pH(t,x,D_x u_\varepsilon,\mu_\varepsilon)) = 0, &(t,x) \in [0,T] \times \Omega^\varepsilon \\
            \mu_\varepsilon(t) = (I,-D_pH(t,\cdot,D_xu_\varepsilon,\mu_\varepsilon)) \# m(t,\cdot), &t \in [0,T] \\
            u_\varepsilon(T,x) = g_\varepsilon(x,m(T)), \hspace{1cm} m_\varepsilon(0,x) = m_0(x), &x \in \Omega^\varepsilon \\
            aD_xu_\varepsilon \cdot \n_\varepsilon = 0, \hspace{1cm} [a^*D_xm_\varepsilon + m_\varepsilon (\Tilde{b}(x) + D_pH(t,x,D_xu_\varepsilon,\mu_\varepsilon))] \cdot \n_\varepsilon = 0, &(t,x) \in [0,T] \times \partial \Omega^\varepsilon
        \end{cases}
        \label{Eq: Approximating System 1}
    \end{equation}
    and
    \begin{equation}
        \begin{cases}
            -\partial_tv_\varepsilon - \sum_{i,j} a_{ij}\partial_{ij}^2 v_\varepsilon + H(t,x,D_xv_\varepsilon,\nu_\varepsilon) = f(t,x,\rho), &(t,x) \in [0,T] \times \Omega^\varepsilon \\
            \partial_t\rho_\varepsilon - \sum_{i,j} \partial_{ij}^2 (a_{ij}\rho_\varepsilon) - \nabla \cdot (\rho_\varepsilon D_pH(t,x,D_xv_\varepsilon,\nu_\varepsilon)) = 0, &(t,x) \in [0,T] \times \Omega^\varepsilon \\
            \nu_\varepsilon(t) = (I,-D_pH(t,\cdot,D_xv_\varepsilon,\nu_\varepsilon)) \# \rho(t,\cdot), &t \in [0,T] \\
            v_\varepsilon(T,x) = g_\varepsilon(x,\rho(T)), \hspace{1cm} \rho_\varepsilon(0,x) = m_0(x), &x \in \Omega^\varepsilon \\
            aD_xv_\varepsilon \cdot \n_\varepsilon = 0, \hspace{1cm} [a^*D_x\rho_\varepsilon + \rho_\varepsilon (\Tilde{b}(x) + D_pH(t,x,D_xv_\varepsilon,\nu_\varepsilon))] \cdot \n_\varepsilon = 0, &(t,x) \in [0,T] \times \partial \Omega^\varepsilon
        \end{cases}
        \label{Eq: Approximating System 2}
    \end{equation}
    respectively. As in section \ref{Sec: Existence 3}, for each $K \subset\subset \Omega$, we have uniform estimates for $u_\varepsilon$ and $v_\varepsilon$ in $C^{\beta/2,1+\beta}([0,T] \times K)$, which implies that every subsequence has a subsequence converging pointwise in $C^{0,1}([0,T] \times \Omega)$. Applying Lemmas \ref{Lem: New Continuity of mu 3} and \ref{Lem: Properties of H} give convergence of $H(t,x,D_xu_{\varphi(\varepsilon)},\mu_{\varphi(\varepsilon)})$, $H(t,x,D_xv_{\varphi(\varepsilon)},\nu_{\varphi(\varepsilon)})$, $D_pH(t,x,D_xu_{\varphi(\varepsilon)},\mu_{\varphi(\varepsilon)})$, and $D_pH(t,x,D_xv_{\varphi(\varepsilon)},\nu_{\varphi(\varepsilon)})$.  By \cite[Corollary 3.9]{porretta2020mean}, it follows that the limits of the subsequences $u_{\varphi(\varepsilon)}, v_{\varphi(\varepsilon)}$ (and hence the sequences themselves) are $u$ and $v$, respectively. Similarly, since we have $D_pH(t,x,D_xu_\varepsilon,\mu_\varepsilon) \rightarrow D_pH(t,x,D_xu,\mu)$ and $D_pH(t,x,D_xv_\varepsilon,\nu_\varepsilon) \rightarrow D_pH(t,x,D_xv,\nu)$ pointwise, \cite[Proposition 4.3]{porretta2020mean} gives that $m_\varepsilon \rightarrow m$ and $\rho_\varepsilon \rightarrow \rho$ in $C^0(0,T;L^1(\Omega))$.
    
    By arguments of the same spirit as those in \cite{kobeissi2022mean}, we get
    $$
    \begin{aligned}
        &\int_{\Omega^\varepsilon} (g_\varepsilon(x,m(T)) - g_\varepsilon(x,\rho(T)))(m_\varepsilon - \rho_\varepsilon) dx + \int_0^T\int_{\Omega^\varepsilon} (f(t,x,m) - f(t,x,\rho))(m_\varepsilon - \rho_\varepsilon)dxdt \\
        &= \int_0^T\int_{\Omega^\varepsilon}  \bigg[m_\varepsilon(H(t,x,D_xu_\varepsilon,\mu_\varepsilon) - H(t,x,D_xv_\varepsilon,\nu_\varepsilon) + D_pH(t,x,D_xu_\varepsilon,\mu_\varepsilon) \cdot (D_xu_\varepsilon - D_xv_\varepsilon)) \\
        &\qquad +\rho_\varepsilon(H(t,x,D_xv_\varepsilon,\nu_\varepsilon) - H(t,x,D_xu_\varepsilon,\mu_\varepsilon) + D_pH(t,x,D_xv_\varepsilon,\nu_\varepsilon) \cdot (D_xv_\varepsilon - D_xu_\varepsilon))\bigg] dxdt \\
        &= \int_0^T\int_{\Omega^\varepsilon}  \bigg[m_\varepsilon(L(t,x,\alpha^{\nu_\varepsilon},\nu_\varepsilon) - L(t,x,\alpha^{\mu_\varepsilon},\mu_\varepsilon) + D_\alpha L(t,x,\alpha^{\nu_\varepsilon},\nu_\varepsilon) \cdot (\alpha^{\mu_\varepsilon} - \alpha^{\nu_\varepsilon})) \\
        &\qquad +\rho_\varepsilon(L(t,x,\alpha^{\mu_\varepsilon},\mu_\varepsilon) - L(t,x,\alpha^{\nu_\varepsilon},\nu_\varepsilon) + D_\alpha L(t,x,\alpha^{\mu_\varepsilon},\mu_\varepsilon) \cdot (\alpha^{\nu_\varepsilon} - \alpha^{\mu_\varepsilon}))\bigg] dxdt \\
        &\leq \int_0^T\int_{\Omega^\varepsilon}  \bigg[m_\varepsilon(L(t,x,\alpha^{\mu_\varepsilon},\nu_\varepsilon) - L(t,x,\alpha^{\mu_\varepsilon},\mu_\varepsilon)) + \rho_\varepsilon(L(t,x,\alpha^{\nu_\varepsilon},\mu_\varepsilon) - L(t,x,\alpha^{\nu_\varepsilon},\nu_\varepsilon))\bigg] dxdt \\
        &= -\int_0^T\int_{\Omega \times \R^n} (L(t,x,\alpha,\mu_\varepsilon) - L(t,x,\alpha,\nu_\varepsilon)) d(\mu_\varepsilon - \nu_\varepsilon) \\
        &\leq 0.
    \end{aligned}
    $$
    By A\ref{A: m0 and f}-\ref{A: g} and the fact that $(m_\varepsilon - \rho_\varepsilon) \rightarrow (m - \rho)$ in $C^0(0,T;L^1(\Omega))$, it follows that
    $$
    \begin{aligned}
        0 &\leq \int_\Omega (g(x,m(T)) - g(x,\rho(T)))(m - \rho) dx + \int_0^T\int_\Omega (f(t,x,m) - f(t,x,\rho))(m - \rho)dxdt \\
        &\leq \int_0^T\int_{\Omega^\varepsilon}  \bigg[m(L(t,x,\alpha^{\mu},\nu) - L(t,x,\alpha^{\mu},\mu)) + \rho_\varepsilon(L(t,x,\alpha^{\nu},\mu) - L(t,x,\alpha^{\nu},\nu))\bigg] dxdt \\
        &\leq 0.
    \end{aligned}
    $$
    In particular, these integrals must vanish. Recall that since $L$ is strictly convex,
    \begin{equation}
        L(t,x,\alpha_1,\mu) - L(t,x,\alpha_2,\mu) + D_\alpha L(t,x,\alpha_1,\mu) \cdot (\alpha_2 - \alpha_1) = 0
    \end{equation}
    if and only if $\alpha_1 = \alpha_2$. Hence, $\alpha^\mu = \alpha^\nu$. By uniqueness of solutions the FP equation, $m = \rho$. By Lemma \ref{Lem: Existence & Bound for mu 2}, this implies $\mu = \nu$. Finally, by uniqueness of solutions to the HJ equation, $u = v$.
\end{proof}

\subsection{Regularity of Solutions}
\label{Sec: Regularity}

In this section, we investigate the regularity of the solutions to our system. First, we state results on the semiconcavity of $u$ and the boundedness of $m$ in $L^\infty(Q)$, which are proved using identical arguments to those used for the analogous results in \cite{porretta2020mean}.

\begin{theorem} \label{Thm: Semiconcavity}
    Assume A\ref{A: Regularity of DpH}-\ref{A: Bound for L and Lx} and A\ref{A: Monotonicity for g,f}-\ref{A: Invariance Condition} hold. Now suppose $\sigma \in W^{2,\infty}(\Omega)$ and
    $$\|f(t,\cdot,m)\|_{W^{2,\infty}(\Omega)} + \|g(\cdot,m)\|_{W^{2,\infty}(\Omega)} \leq C_0$$
    for all $t \in [0,T]$ and $m \in L^1(\Omega)$. Finally, assume that for $\nu \in \mathfrak{M}_{\infty,R}(\Omega \times \R^n)$, there exist constants $C_1,C_2$, depending only on $\Lambda_{q'}(\nu)$, such that
    $$
    \begin{aligned}
        &H(t,x,p,\nu) + H(t,y,q,\nu) - 2H\left(t,z,\frac{p+q}{2},\nu\right) \\
        &\geq -C_1(|x-z|^2 + |y-z|^2 + |x+y-2z|^2)(1 + |p+q|) - C_2|x-y||p-q|
    \end{aligned}
    $$
    for any $x,y,z \in \Omega$, $p,q \in \R^n$, and $t \in (0,T)$. If $(u,m,\mu)$ is a solution to \eqref{Eq: MFGC 2}, then $u(t,\cdot)$ is semiconcave for all $t \in (0,T)$, with a semiconcavity constant bounded uniformly in $t$. In particular,
    $$D_{xx}^2u(t,\cdot) \leq M$$
    for all $t \in (0,T)$.
\end{theorem}

\begin{theorem} \label{Thm: Regularity of m}
    Assume the conditions of Theorem \ref{Thm: Semiconcavity} hold. Now suppose that there exists $\delta_0 > 0$ such that
    $$\left(\Tilde{b}(x) + D_pH(t,x,p,\mu)\right) \cdot D_xd(x) \leq 0$$
    for all $x \in \Gamma_{\delta_0}$ and $p \in \R^n$. Then if $(u,m,\mu)$ is a solution to \eqref{Eq: MFGC 2}, we have $m \in L^\infty(Q)$ with
    $$\|m\|_{L^\infty(Q)} \leq C.$$
\end{theorem}

Finally, we prove that under a few additional assumptions, the value function $u$ will satisfy the HJ equation in a classical sense on the interior of our domain.

\begin{theorem} \label{Thm: Regularity of solutions}
    Assume A\ref{A: Regularity of DpH}-\ref{A: Bound for L and Lx} and A\ref{A: Monotonicity for g,f}-\ref{A: Invariance Condition} hold. Now assume that for $\rho \in L^1(\Omega)$ with $\|\rho\|_{L^1(\Omega)} \leq 1$, we have $\|g_\varepsilon(\cdot,m)\|_{C^{2+\beta}(\Omega)} \leq C_0$. Finally, assume there exist $L_1 \in (0,1)$ and $\xi_1: \R^2 \rightarrow [0,\infty)$ continuous such that for $R > 0$, $p \in \R^n$, and $\mu_1,\mu_2 \in \mathfrak{M}_{\infty,R}(\Omega \times \R^n)$ we have:
    \begin{itemize}
        \item If $\mu_i = (I,\alpha^{\mu_i}) \# \rho$, then
        $$|D_pH(t,x,p,\mu_1) - D_pH(t,x,p,\mu_2)| \leq L_1\|\alpha^{\mu_1} - \alpha^{\mu_2}\|_{L^{q_0}(\rho)}.$$
        \item If $\mu_i = (I,\alpha) \# \rho_i$, then
        $$|D_pH(t,x,p,\mu_1) - D_pH(t,x,p,\mu_2)| \leq d^*(\rho_1,\rho_2)^\beta\xi_1(\abs{p},\|\alpha\|_\infty).$$
    \end{itemize}
    Then the unique weak solution $(u,m,\mu)$ to \eqref{Eq: MFGC 2} is actually a classical solution (in the sense that $u$ is a classical solution to the HJ equation in $[0,T] \times \Omega$).
\end{theorem}

\begin{proof}
    As in the proof of Theorem \ref{Thm: Existence}, let $(u_\varepsilon,m_\varepsilon,\mu_\varepsilon)$ be a sequence of solutions to \eqref{Eq: Approximating MFGC} converging to $(u,m,\mu)$. However, in this case, the $(u_\varepsilon,m_\varepsilon,\mu_\varepsilon)$ are classical solutions on $[0,T] \times \overline{\Omega}_\varepsilon$. Combining the Lipschitz estimate from Lemma \ref{Lem: Gradient Bound} with the arguments from Sections \ref{Sec: Bootstrapping 1} and \ref{Sec: Bootstrapping 2}, for each $K \subset\subset \Omega$, we get bounds for $u_\varepsilon$ in $C^{1+\alpha/2,2+\alpha}([0,T] \times K)$ uniformly in $\varepsilon$. By the Arzela-Ascoli theorem and a diagonal argument, we get a subsequence that converges in $C^{1,2}([0,T] \times K)$ for all $K \subset\subset \Omega$, and hence $u_\varepsilon$ and its derivatives converge pointwise in $[0,T] \times \Omega$. By the uniqueness of the limit, we get that $u \in C^{1,2}([0,T] \times K)$ for all $K \subset\subset \Omega$. Furthermore, by A\ref{A: g} and Lemma \ref{Lem: Properties of H}, passing to the limit gives that $u$ satisfies the HJ equation on $[0,T] \times \Omega$.
\end{proof}

\bibliographystyle{amsplain}
\bibliography{refs}

\end{document}